\documentclass[a4paper,12pt]{article}
\usepackage{latexsym,bm}
\usepackage{amsmath,amsfonts,amssymb}

\newcommand{\BS}{\mathfrak S}\renewcommand{\l}{\ell}
\newcommand{\Z}{\mathbb Z} \renewcommand{\H}{\mathcal{H}}
 
 \newcommand{\lam}{\lambda}
  
  \DeclareMathOperator{\Ind}{Ind}

 \DeclareMathOperator{\res}{res}
\DeclareMathOperator{\mmod}{mod} \DeclareMathOperator{\id}{id}
\newcommand{\ts}{\widetilde{S}} \DeclareMathOperator{\Aut}{Aut}
\DeclareMathOperator{\GL}{GL}
\newcommand{\eps}{\varepsilon} \DeclareMathOperator{\rad}{rad}
\newcommand{\veps}{\varepsilon}
 \newcommand{\td}{\widetilde{D}}
\DeclareMathOperator{\ch}{char} \newcommand{\A}{\mathcal A}
\newcommand{\ulam}{{\lam}} \newcommand{\sig}{\sigma}
\newcommand{\umu}{{\mu}} 
\DeclareMathOperator{\aff}{aff} \DeclareMathOperator{\HH}{h}
\DeclareMathOperator{\bQ}{Q}
\newcommand{\bbQ}{\overrightarrow\bQ} \DeclareMathOperator{\rres}{res}
\newcommand{\kgl}{\widehat{\mathfrak{gl}}} \DeclareMathOperator{\Irr}{Irr}
\newcommand{\ksl}{\widehat{\mathfrak{sl}}}
\newcommand{\qed}{\hfill\text{$\square$}}

\newtheorem{prop}{Proposition}
\newtheorem{thm}{Theorem}\newtheorem{cor}{Corollary}
\newtheorem{lem}{Lemma}\newtheorem{dfn}{Definition}
\numberwithin{equation}{section} \numberwithin{prop}{section}
\numberwithin{thm}{section} \numberwithin{lem}{section}
\numberwithin{dfn}{section} \numberwithin{cor}{section}

\title{Crystal bases and simple modules for Hecke algebras of type
$G(p,p,n)$
\thanks{Research supported by National Natural Science Foundation of
China (Project 10401005) and by Program for New Century Excellent
Talents in University and partly by the URF of Victoria University
of Wellington. } \footnotetext{{\it 2000 Mathematics Subject
Classification.} Primary 20C08, 20C20, 17B37.} \footnotetext{{\it
Key words and phrases.} Cyclotomic Hecke algebra, Fock space,
crystal basis, Kleshchev multipartition, Lakshmibai--Seshadri
path.}}
\author{Jun Hu \\[5pt]
Department of Applied Mathematics\\
Beijing Institute of Technology\\
Beijing, 100081, P.R. China\null\\[1.5pt]
E-mail: junhu303@yahoo.com.cn}
\date{}

\begin{document}
\maketitle
\begin{abstract}

We apply the crystal basis theory for Fock spaces over quantum
affine algebras to the modular representations of the cyclotomic
Hecke algebras of type $G(p,p,n)$. This yields a classification of
simple modules over these cyclotomic Hecke algebras in the
non-separated case, generalizing our previous work [J. Hu, J.
Algebra 267 (2003) 7-20]. The separated case was completed in [J.
Hu, J. Algebra 274 (2004) 446--490]. Furthermore, we use
Naito--Sagaki's work [S. Naito \& D. Sagaki, J. Algebra 251 (2002)
461--474] on Lakshmibai--Seshadri paths fixed by diagram
automorphisms to derive explicit formulas for the number of simple
modules over these Hecke algebras. These formulas generalize earlier
results of [M. Geck, Represent. Theory 4 (2000) 370-397] on the
Hecke algebras of type $D_n$ (i.e., of type $G(2,2,n)$).
\end{abstract}

\section{Introduction}

The theory of crystal (or canonical) bases is one of the most
significant advances in Lie theory over the past two decades. It was
discovered and developed by M. Kashiwara (\cite{Kas}) and G. Lusztig
(\cite{Lu}) around 1990. Since then remarkable applications to some
classical problems in representation theory have been found. One
typical example is the well-known Lascoux--Leclerc--Thibon's
Conjecture (\cite{LLT}), which asserts that, the decomposition
numbers of the Iwahori--Hecke algebra associated to symmetric group
at a primitive $e$-th root of unity in $\mathbb{C}$ (the complex
number field) can be obtained from the evaluation at $1$ of the
coefficient polynomials of natural bases appeared in the expansion
of global crystal bases (i.e., canonical bases) of some level one
Fock spaces over the quantum affine algebra $U_q(\ksl_e)$. This
conjecture has been proved by S. Ariki (\cite{A1}), who generalized
it to the case of the cyclotomic Hecke algebras of type $G(r,1,n)$.
A similar conjecture (\cite{LT}), which relates the decomposition
numbers of the $q$-Schur algebra with $q$ specialized to a primitive
$e$-th root of unity in $\mathbb{C}$ to global crystal bases of Fock
space as $U_q(\kgl_e)$-module, has been proved by
Varagnolo--Vasserot (\cite{VV}). For further example, see the work
of Brundan--Kleshchev (\cite{BK}), where the theory of crystal bases
of type $A_{2l}^{(2)}$ was applied to the modular representations of
Hecke--Clifford superalgebras as well as of double covers of
symmetric groups.
\smallskip

This paper provides a new application of the theory of crystal (or
canonical) bases to modular representation theory. Precisely, we
apply the crystal basis theory for Fock spaces of higher level over
the quantum affine algebra of type $A_{l}^{(1)}$ to the modular
representations of the cyclotomic Hecke algebra $\H(p,p,n)$ of type
$G(p,p,n)$ in the non-separated case (see Definition \ref{sep}). The
separated case has been completely solved in our previous work
\cite{Hu3}. We explicitly describe (in terms of combinatorics over
certain Kleshchev's good lattices) which irreducible representation
of the Ariki--Koike algebra $\H(p,n)$ remains irreducible when
restricted to $\H(p,p,n)$. This yields a classification of simple
modules over the cyclotomic Hecke algebra $\H(p,p,n)$ in the
non-separated case, generalizing our previous work \cite{Hu2} on the
Hecke algebras of type $D_n$. Then we go further in the remaining
part of this paper. We use Naito--Sagaki's work
(\cite{NS1},\cite{NS2}) on Lakshmibai--Seshadri paths fixed by
diagram automorphisms to derive explicit formula for the number of
simple modules over the cyclotomic Hecke algebra $\H(p,p,n)$. Our
formula generalizes earlier result of Geck \cite{Ge} on the Hecke
algebra of type $D_n$ (i.e., of type $G(2,2,n)$). Note that our
approach even in that special case is quite different, because it is
based on Ariki's celebrated theorem (\cite{A1}) on a generalization
of Lascoux--Leclerc--Thibon's Conjecture as well as Naito--Sagaki's
work (\cite{NS1},\cite{NS2}) on Lakshmibai--Seshadri paths, while
Geck's method in \cite{Ge} depends on explicit information on
character tables and Kazhdan--Lusztig theory for Iwahori--Hecke
algebras associated to finite Weyl groups---not presently available
in our general $G(p,p,n)$ cases. As a byproduct, we get a remarkable
bijection between two sets of Kleshchev multipartitions. Our
explicit formulas strongly indicate that there are some new intimate
connections between the representation of $\H(p,p,n)$ at roots of
unity and the representations of various Ariki--Koike algebras of
smaller sizes at various roots of unity. Although we will not
discuss these matters in the present paper, we remark that it seems
very likely the decomposition matrix of the latter can be naturally
embedded as a submatrix of the decomposition matrix of the former.
\smallskip

The paper is organized as follows. Section 2 collects some basic
known results about Ariki--Koike algebras (i.e., the cyclotomic
Hecke algebras of type $G(r,1,n)$). These include
Dipper--James--Mathas's work on the structure and representation
theory of Ariki--Koike algebras as well as Dipper--Mathas's Morita
equivalence results. The notion of Kleshchev multipartition as well
as Ariki's remarkable result (Theorem \ref{Ariki}) are also
introduced there. In Section 3, we first recall our previous work on
modular representations of Hecke algebras of type $D_n$ and of type
$G(p,p,n)$. Then we give the first two main results (Theorem
\ref{main1} and Theorem \ref{main2}) in this paper. Theorem
\ref{main3} shows that these two main results are valid over any
field $K$ which contains primitive $p$-th root of unity and over
which $\H(p,p,n)$ is split. Our Theorem \ref{main1} is a direct
generalization of \cite[(1.5)]{Hu2}. A sketch of the proof
(following the streamline of the proof of \cite[(1.5)]{Hu2}) is
presented in Section 4. The proof of Theorem \ref{main2} is given in
Section 5. Our main tools used there are Dipper--Mathas's Morita
equivalence results for Ariki--Koike algebras and their connections
with type $A$ affine Hecke algebras. In Section 6 we give the second
two main results (Theorem \ref{mainthm3} and Theorem \ref{mainthm4})
in this paper, which yield explicit formula for the number of simple
modules over the cyclotomic Hecke algebra $\H(p,p,n)$ in the
non-separated case. Note that in the separated case such a formula
can be easily written down by using the results \cite[(5.7)]{Hu3}.
The proof uses our first two main results (Theorem \ref{main1} and
Theorem \ref{main2}) as well as Naito--Sagaki's work
(\cite{NS1},\cite{NS2}) on Lakshmibai--Seshadri paths fixed by
diagram automorphisms. As a byproduct, we get (Corollary
\ref{cor609}) a remarkable bijection between two sets of Kleshchev
multipartitions, which seems of independent interest.
\smallskip

The present paper is an expanded version of an earlier preprint
(cited in Ariki's book \cite[{[}cyclohecke12{]}]{A3}) completed in
the February of 2002. That preprint already contains Theorem
\ref{main1} and Theorem \ref{main2}, which are generalizations of
\cite[(1.5)]{Hu2}. Part of the remaining work was done during the
author's visit of RIMS in 2004. After this expanded version was
completed and the main results were announced, N. Jacon informed us
that a result similar to Theorem \ref{main1} in the context of FLOTW
partitions (see \cite[Definition 2.2]{Ja}) was also obtained in his
Ph.D. Thesis \cite{Ja0} in 2004, and we are informed the existence
of a preprint \cite{GJ} of Genet and Jacon on the modular
representation of the Hecke algebra of type $G(r,p,n)$ (where
$p|r$). Our main results Theorem \ref{mainthm3} and Theorem
\ref{mainthm4} are not related to any results in \cite{Ja0} and
\cite{GJ}. Both the paper \cite{GJ} and the present paper use
Ariki's results on Fock spaces, crystal graphs as well as Clifford
theory. But \cite{GJ} uses a different version of Fock space and
hence a different parameterization of simple modules over
Ariki--Koike algebras. The relationships between the
parameterization results given in \cite{GJ} and our parameterization
results given in Theorem \ref{main1}, Theorem \ref{main2} and in
\cite[Theorem 4.9]{Hu3} are explained in Remark 3.12. \medskip

The author thanks Professor Susumu Ariki for some stimulating
discussion, especially for informing him about the work of
Naito--Sagaki. The author also thanks Professor Masaki Kashiwara for
explaining a result about crystal bases. The main results of this
paper were announced at the ``International Conference on
Representation Theory, III" (Chengdu, August 2004). The author would
like to thank Professor Nanhua Xi and Professor Jie Du for their
helpful comments.  The author also would like to thank the referee
for many helpful suggestions.
\medskip\bigskip

\section{Preliminaries}

Let $r$, $p$, $d$ and $n$ be positive integers such that $pd=r$.
The complex reflection group $G(r,p,n)$ is the group consisting of
$n$ by $n$ permutation matrices with the properties that the
entries are either $0$ or $r$-th roots of unity in $\mathbb{C}$,
and the $d$-th power of the product of the non-zero entries of
each matrix is $1$. The order of $G(r,p,n)$ is $dr^{n-1}n!$, and
$G(r,p,n)$ is a normal subgroup of $G(r,1,n)$ of index $p$.

Cyclotomic Hecke algebras associated to complex reflection groups
were introduced in the work of Brou\'e--Malle (\cite{BM}) and in the
work of Ariki--Koike (\cite{AK}). These algebras are deformations of
the group rings of the complex reflection groups. We recall their
definitions. Let $K$ be a field and let $q, Q_1,\cdots,Q_r$ be
elements of $K$ with $q$ invertible. Let
$\H_K(r,n)=\H_{q,Q_1,\cdots,Q_r}(r,n)$ be the unital $K$-algebra
with generators $T_0,T_1,\cdots,T_{n-1}$ and relations
$$\begin{aligned}
&(T_0-Q_1)\cdots (T_0-Q_r)=0,\\
&T_0T_1T_0T_1=T_1T_0T_1T_0,\\
&(T_i+1)(T_i-q)=0,\quad\text{for $1\leq i\leq n-1$,}\\
&T_iT_{i+1}T_i=T_{i+1}T_{i}T_{i+1},\quad\text{for $1\leq i\leq n-2$,}\\
&T_iT_j=T_jT_i,\quad\text{for $0\leq i<j-1\leq n-2$.}\end{aligned}
$$
This algebra is called {\it Ariki--Koike algebra} or the cyclotomic
Hecke algebra of type $G(r,1,n)$. Whenever the parameter $q$ is
clear from the context, we shall say (for simplicity) that the
algebra $\H_K(r,n)$ is with parameters set $\{Q_1,\cdots,Q_r\}$.
This algebra contains the Hecke algebras of type $A$ and type $B$ as
special cases. It can be defined over $\Z[v,v^{-1},v_1,\cdots,v_r]$,
where $v, v_1,\cdots, v_r$ are all indeterminates. Upon setting
$v=1$ and $v_i=(\sqrt[r]{1})^{i-1}$ for each $i$, where
$\sqrt[r]{1}$ denotes a primitive $r$-th root of unity in
$\mathbb{C}$, one obtains the group algebra for the complex
reflection group $G(r,1,n)\cong\Z_{r}\wr\BS_{n}$.
\smallskip

Suppose that $K$ contains a primitive $p$-th root of unity $\eps$.
Let $x_1,\cdots,x_d$ be invertible elements in $K$ with
$x_{i}^{1/p}\in K$ for each $i$. We consider the Hecke algebra
$\H_{K}(r,n)$ with parameters
$$ q,\,\, x_{i}^{1/p}\eps^{j},\quad
i=1,2,\cdots,d,\,\,j=0,1,\cdots,p-1.
$$
Then the first defining relation for $\H_{K}(r,n)$ becomes $$
(T_{0}^{p}-x_1)(T_{0}^{p}-x_2)\cdots (T_{0}^{p}-x_d)=0. $$ Let
$\H_{K}(r,p,n)=\H_{q,x_{1},\cdots,x_{d}}(r,p,n)$ be the subalgebra
of $\H_{K}(r,n)$ generated by the elements $$ T_0^{p},\quad
T_{u}:{=}T_{0}^{-1}T_{1}T_{0}, \,\,\,
T_{1},\,T_2,\,\cdots,\,T_{n-1}. $$ Then it is a $q$-analogue of the
group algebra for the complex reflection group $G(r,p,n)$. This
algebra is called the cyclotomic Hecke algebra of type $G(r,p,n)$.
It is known that in this case (by \cite{MaM}) $\H_{K}(r,p,n)$ is a
symmetric algebra over $K$. For simplicity, we shall often write
$\H(r,p,n), \H(r,n)$ instead of $\H_{K}(r,p,n), \H_{K}(r,n)$.

Our main interest in this paper will be the algebra $\H(r,p,n)$ in
the special case when $p=r$, that is, the cyclotomic Hecke algebra
of type $G(p,p,n)$. The special case $\H(2,2,n)$ is just the
Iwahori--Hecke algebra of type $D_n$. For convenience, we shall use
a normalized version of $\H(p,p,n)$ which is defined as follows. Let
$p, n\in\mathbb{N}$. Let $K$ be a field and $q$ be an invertible
element in $K$. Throughout this paper, we assume that $K$ contains
primitive $p$-th roots of unity. In particular, $\ch K$ is coprime
to $p$. We define the algebra $\H_q(p,n)$ to be the associative
unital $K$-algebra with generators $T_0, T_1,\cdots, T_{n-1}$
subject to the following relations
$$\begin{aligned}
&T_{0}^{p}-1=0,\\
&T_0T_1T_0T_1=T_1T_0T_1T_0,\\
&(T_i+1)(T_i-q)=0,\quad\text{for $1\leq i\leq n-1$,}\\
&T_iT_{i+1}T_i=T_{i+1}T_{i}T_{i+1},\quad\text{for $1\leq i\leq n-2$,}\\
&T_iT_j=T_jT_i,\quad\text{for $0\leq i<j-1\leq n-2$.}\end{aligned}
$$
Note that if we fix a choice of a primitive $p$-th root of unity
$\eps$ in $K$, then the first relation becomes
$$(T_0-1)(T_0-\eps)\cdots(T_{0}-\eps^{p-1})=0,$$
and the algebra $\H_q(p,n)$ is just the Ariki--Koike algebra or the
cyclotomic Hecke algebra of type $G(p,1,n)$ with parameters
$q,1,\eps,\cdots,\eps^{p-1}$. However, the algebra $\H_q(p,n)$
itself does not depend on the choice of primitive $p$-th root of
unity. Now, let $\H_{q}(p,p,n)$ be the subalgebra of $\H_q(p,n)$
generated by the elements
$$
T_{u}:=T_{0}^{-1}T_{1}T_{0},\,\,\, T_{1},\,T_2,\,\cdots,\,T_{n-1}.
$$
This algebra, called cyclotomic Hecke algebra of type $G(p,p,n)$,
will be the main subject of this paper.\smallskip

It is well-known that $\H_{q}(p,n)$ is a free
$\H_{q}(p,p,n)$-module with basis $\{1,
T_{0},\cdots,T_{0}^{p-1}\}$. As a $\H_{q}(p,p,n)$-module,
$\H_{q}(p,n)$ is in fact isomorphic to a direct sum of $p$ copies
of regular $\H_{q}(p,p,n)$-modules. Let $\tau$ be the $K$-algebra
automorphism of $\H_{q}(p,n)$ which is defined on generators by
$\tau(T_1)=T_0^{-1}T_1T_0, \tau(T_i)=T_i$, for any $i\neq 1$. Let
$\sig$ be the nontrivial $K$-algebra automorphism of $\H_{q}(p,n)$
which is defined on generators by $\sig(T_0)=\eps T_0,
\sig(T_i)=T_i$, for any $1\leq i\leq n-1$. By \cite[(1.4)]{Hu3},
$\tau(\H_{q}(p,p,n))=\H_{q}(p,p,n)$ and clearly
$\sig\downarrow_{\H_{q}(p,p,n)}=\id$. Moreover, the set of
$K$-subspaces $\bigl\{T_0^{i}\H_{q}(p,p,n)\bigr\}_{i=0}^{p-1}$ of
$\H_{q}(p,n)$ forms a $\Z/p\Z$-graded Clifford system in
$\H_{q}(p,n)$ in the sense of \cite[(11.12)]{CR}.\smallskip

Our approach to the modular representations of the algebra
$\H_q(p,p,n)$ is to consider the restriction of the representations
of $\H_q(p,n)$. To this end, we have to first recall some known
results about $\H_q(p,n)$. The structure and representation theory
for Ariki--Koike algebras with arbitrary parameters $q,Q_1,$
$\cdots,Q_p$ have been well studied in \cite{DJM}, where it was
shown that these algebras are cellular in the sense of \cite{GL} if
$q$ is invertible. To state their results, we need some
combinatorics. Recall that a partition of $n$ is a non-increasing
sequence of positive integers $\lam=(\lam_1,\cdots,\lam_r)$ such
that $|\lam|:=\sum_{i=1}^{r}\lam_i=n$, while a $p$-multipartition of
$n$ is a $p$-tuple of partitions
$\ulam=(\lam^{(1)},\cdots,\lam^{(p)})$ such that
$|\lam|:=\sum_{i=1}^p|\lam^{(i)}|=n$. For any two
$p$-multipartitions $\ulam, \umu$ of $n$, we define
$\ulam\trianglerighteq\umu$ if
$$ \sum_{j=1}^{i-1}|\lam^{(j)}|+\sum_{k=1}^{m}\lam^{(i)}_k \geq
\sum_{j=1}^{i-1}|\mu^{(j)}|+\sum_{k=1}^{m}\mu^{(i)}_k,\,\,\,\forall\,1\leq
i\leq p,\,m\geq 1.
$$
We shall state the results of \cite{DJM} in our special case of the
Ariki--Koike algebra $\H_{q}(p,n)$. First, we fix a primitive $p$-th
root of unity $\eps$ in $K$. Let $
\bQ:=\{1,\eps,\eps^2,\cdots,\eps^{p-1}\}$. Let
$\bbQ=(Q_1,\cdots,Q_p)$ be an ordered $p$-tuple which is obtained by
fixing an order on $\bQ$. We regard $\H_{q}(p,n)$ as the
Ariki--Koike algebra with parameters $q,1,\eps,\cdots,$
$\eps^{p-1}$. Then for each $p$-multipartitions $\ulam$ of $n$,
there is a Specht module, denoted by $\ts^{\ulam}_{\bbQ}$, and there
is a naturally defined bilinear form $\langle,\rangle$ on
$\ts^{\ulam}_{\bbQ}$. Let
$\td^{\lam}_{\bbQ}=\ts^{\ulam}_{\bbQ}/\rad\langle,\rangle$. Note
here the subscript $\bbQ$ (instead of $\bQ$) is used to emphasis
that the structure of the module $\ts^{\ulam}_{\bbQ}$ (resp.,
$\td^{\lam}_{\bbQ}$) does depend on the fixed order on the set
$\{1,\eps,\cdots,\eps^{p-1}\}$ of parameters.

\begin{dfn} \label{def21} Let $v$ be an indeterminate over $\Z$. Let $\epsilon$ be a
primitive $p$-th root of unity in $\mathbb{C}$. We define
$\mathcal{Z}=\Z[\epsilon][v,v^{-1}]$ and
$$ f_{p,n}(v,\epsilon)=\prod_{1\leq i<j\leq
p}\prod_{-n<k<n}\bigl(\epsilon^{i-1}v^k-\epsilon^{j-1}\bigr)\in
\mathcal{Z}.
$$
\end{dfn}

Note that $v^{p}-1=\prod_{k|p}\Phi_k(v)$, where $\Phi_k(v)$ is the
$k$-th cyclotomic polynomial over $\Z$. It follows easily that for
any $\Z[v,v^{-1}]$-algebra $K$ which contains a primitive $p$-th
root of unity $\eps$, the natural homomorphism
$\Z[v,v^{-1}]\rightarrow K$ can be uniquely extended to a
homomorphism from $\mathcal{Z}$ to $K$ by mapping $\epsilon$ to
$\eps$. In other words, we can always specialize $\epsilon$ to any
primitive $p$-th root of unity.

\addtocounter{thm}{1}
\begin{thm} {\rm (\cite{DJM}, \cite{A4})}\label{thm22} With the above notations, we have that \smallskip

1) the set $\bigl\{\td^{\ulam}_{\bbQ}\bigm|\text{$\ulam$ is a
$p$-multipartition of $n$ and $\td_{\bbQ}^{\ulam}\neq 0$}\bigr\}$
forms a complete set of pairwise non-isomorphic simple
$\H_{q}(p,n)$-modules;\smallskip

2) if $\td_{\bbQ}^{\umu}\neq 0$ is a composition factor of
$\ts_{\bbQ}^{\ulam}$ then $\ulam\trianglerighteq\umu$, and every
composition factor of $\ts_{\bbQ}^{\ulam}$ is isomorphic to some
$\td_{\bbQ}^{\umu}$ with $\ulam\trianglerighteq\umu$; if
$\td_{\bbQ}^{\ulam}\neq 0$ then the composition multiplicity of
$\td_{\bbQ}^{\ulam}$ in $\ts_{\bbQ}^{\ulam}$ is $1$;
\smallskip

3) $\H_q(p,n)$ is semisimple if and only if
$$\biggl(\prod_{i=1}^{n}\bigl(1+q+q^2+\cdots+q^{i-1}\bigr)\biggr)f_{p,n}(q,\eps)\neq
0,$$ in $K$. In that case, $\ts_{\bbQ}^{\lam}=\td_{\bbQ}^{\ulam}$
for each $p$-multipartition $\ulam$ of $n$.
\end{thm}
 Note that whether $f_{p,n}(q,\eps)$ is nonzero in $K$ or not is
independent of the choice of the primitive $p$-th root of unity
$\eps$ in $K$.\medskip

It remains to determine when $\td_{\bbQ}^{\ulam}\neq 0$. This was
solved in \cite{DJ} in the case of type $A$. In general case it was
solved by the work of Dipper--Mathas \cite{DM} and the work of Ariki
\cite{A2}. We first recall Dipper--Mathas's result in \cite{DM}. Two
parameters $Q_i, Q_j$ are said to be in the same $q$-orbit, if
$Q_i=q^kQ_j$ for some $k\in\Z$. Now we suppose that
$\bQ=\bQ_1\sqcup\bQ_2\sqcup\cdots\sqcup\bQ_{\kappa}$ (disjoint
union) such that $Q_i, Q_j$ are in the same $q$-orbit only if $Q_i,
Q_j\in\bQ_c$ for some integer $c$ with $1\leq c\leq\kappa$. Let
$p_i=|\bQ_i|$ for each integer $i$ with $1\leq i\leq\kappa$.

\begin{thm} {\rm (\cite[Theorem 1.1]{DM})}\label{thm23} With the above notations, the algebra \smallskip
$\H_q(p,n)$ is Morita equivalent to the algebra $$
\bigoplus_{\substack{n_1,\cdots,n_{\kappa}\geq 0\\
n_1+\cdots+n_{\kappa}=n}}\H_{q,\bQ_1}(p_1,n_1)\otimes\cdots\otimes\H_{q,\bQ_{\kappa}}(p_{\kappa},n_{\kappa}),
$$
where each $\H_{q,\bQ_i}(p_i,n_i)$ denotes the Ariki--Koike algebra
of size $n_i$ and with parameters set $\bQ_i$. Moreover, if we fix
an order on each $\bQ_i$ to get an ordered tuple $\bbQ_i$ and
suppose that $\bbQ=\bigl(\bbQ_1,\bbQ_2,\cdots,\bbQ_{\kappa}\bigr)$
(concatenation of ordered tuples), then the above Morita equivalence
sends $\td_{\bbQ}^{\lam}$ to
$\td_{\bbQ_1}^{\lam^{[1]}}\otimes\cdots\otimes\td_{\bbQ_{\kappa}}^{\lam^{[\kappa]}}$,
where $$
\lam^{[i]}:=\bigl(\lam^{(\sum_{j=1}^{i-1}p_j+1)},\lam^{(\sum_{j=1}^{i-1}p_j+2)},\cdots,\lam^{(\sum_{j=1}^{i}p_j)}\bigr),
\,\,|\lam^{[i]}|=n_i,\,\,\forall\,1\leq i\leq\kappa, $$
$\td_{\bbQ_i}^{\lam^{[i]}}$ denotes the quotient module (see
definition above (\ref{def21})) of the Specht module
$\ts_{\bbQ_i}^{\lam^{[i]}}$ over $\H_{q,\bQ_i}(p_i,n_i)$. In
particular, $\td_{\bbQ}^{\lam}\neq 0$ if and only if
$\td_{\bbQ_i}^{\lam^{[i]}}\neq 0$ for any integer $i$ with $1\leq
i\leq\kappa$.
\end{thm}

\addtocounter{lem}{3}
\begin{lem} \label{lm24} Let $\H_{q,\bQ}(p,n)$ be the Ariki--Koike algebra with
parameters $q,Q_1,$ $\cdots,Q_p$ and defined over $K$. Let
$\bQ=\{Q_1,\cdots,Q_p\}$. Let $0\neq a\in K$. Let
$a\!\bQ=\{aQ_1,\cdots,aQ_p\}$. Let $\sigma_a$ be the isomorphism
from $\H_{q,a\!\bQ}(p,n)$ onto $\H_{q,\bQ}(p,n)$ which is defined on
generators by $\sigma_a(T_0)=aT_0$ and $\sigma_a(T_i)=T_i$ for $
i=1,2,\cdots,n-1$. Let $\bbQ$ be an ordered $p$-tuple which is
obtained by fixing an order on $\bQ$. Then for each
$p$-multipartition $\lam$ of $n$, there are
$\H_{q,a\!\bQ}(p,n)$-module isomorphisms
$$\Bigl(\ts_{\bbQ}^{\lam}\Bigr)^{\sigma_a}\cong\ts_{\overrightarrow{a\!\bQ}}^{\lam},\,\,\,\,
\Bigl(\td_{\bbQ}^{\lam}\Bigr)^{\sigma_a}\cong\td_{\overrightarrow{a\!\bQ}}^{\lam},$$
where $\overrightarrow{a\!\bQ}$ denotes the ordered $p$-tuple
which is obtained from $\bbQ$ by multiplying $a$ on each
component. In particular, $\td_{\bbQ}^{\lam}\neq 0$ if and only if
$\td_{\overrightarrow{a\!\bQ}}^{\lam}\neq 0$.
\end{lem}

\noindent {Proof:} \,This follows directly from the definitions
and constructions of the modules $\ts_{\bbQ}^{\lam},
\td_{\bbQ}^{\lam}$.\hfill\qed
\medskip

Theorem \ref{thm23} and Lemma \ref{lm24} reduce the problem on
determining when $\td_{\bbQ}^{\ulam}\neq 0$ to the case where all
the parameters in $\bQ$ are some powers of $q$. In that case, the
problem was solved by Ariki \cite{A2}. From now on until the end of
this section, we assume that all the parameters in $\bQ$ are some
powers of $q$. In particular, in our $\H_q(p,n)$ case, $q$ must be a
root of unity. To state the result of Ariki, we have to recall the
definition of Kleshchev multipartition (see \cite{AM}). For any
$p$-multipartition $\lam$, the Young diagram of $\lam$ is the set
$$ [\lam]=\Bigl\{(a,b,c)\Bigm|\text{$1\leq c\leq p$ and $1\leq
b\leq\lam_{a}^{(c)}$}\Bigr\}. $$ The elements of $[\lam]$ are
nodes of $\lam$. Given any two nodes $\gamma=(a,b,c),
\gamma'=(a',b',c')$ of $\lam$, say that $\gamma$ is {\it below}
$\gamma'$, or $\gamma'$ is {\it above} $\gamma$, if either $c>c'$
or $c=c'$ and $a>a'$. The {\it residue} of $\gamma=(a,b,c)$ is
defined to be \addtocounter{equation}{4}
\begin{equation} \label{res}\text{$\rres(\gamma):= m+e\Z\in\Z/e\Z$,
\quad if $q=\sqrt[e]{1}$ and $q^m= q^{b-a}Q_c$,}\end{equation} and
we say that $\gamma$ is a $\rres(\gamma)$-node. A {\it removable}
node is a node of the boundary of $[\lam]$ which can be removed,
while an {\it addable} node is a concave corner on the rim of
$[\lam]$ where a node can be added. If
$\mu=(\mu^{(1)},\cdots,\mu^{(p)})$ is a $p$-multipartition of
$n+1$ with $[\mu]=[\lam]\cup\bigl\{\gamma\bigr\}$ for some
removable node $\gamma$ of $\mu$, we write $\lam\rightarrow\mu$.
If in addition $\rres(\gamma)=x$, we also write that
$\lam\overset{x}{\rightarrow}\mu$. For example, suppose $n=10,
p=4, q=\sqrt[8]{1}$ and $\eps=q^2=\sqrt[4]{1}$. The nodes of
$\lam=((2,1),(1^2),(1^3),(2))$ have the following residues $$
\lam=\biggl(\left(\begin{matrix} 0& 1\\
7
\end{matrix}
\right),\left(\begin{matrix} 2 \\
1
\end{matrix}
\right),\left(\begin{matrix} 4 \\
3\\
2
\end{matrix}
\right),\left(\begin{matrix} 6& 7\end{matrix} \right)\biggr).
$$
It has five removable nodes. Fix a residue $x$ and consider the
sequence of removable and addable $x$-nodes obtained by reading the
boundary of $\lam$ from the bottom up. In the above example, we
consider residue $x=1$, then we get a sequence ``ARR", where each
``A'' corresponds to an addable $x$-node and each ``R'' corresponds
to a removable $x$-node. Given such a sequence of letters ``A,R" ,
we remove all occurrences of the string ``AR'' and keep on doing
this until no such string ``AR'' is left. The ``R''s that still
remain are the {\it normal} $x$-nodes of $\lam$ and the highest of
these is the {\it good} $x$-node. In the above example, there is
only one normal $1$-node, which is a good $1$-node. If $\gamma$ is a
good $x$-node of $\mu$ and $\lam$ is the multipartition such that
$[\mu]=[\lam]\cup\gamma$, we write
$\lam\overset{x}{\twoheadrightarrow}\mu$. For each integer $n\geq
0$, let $\mathcal{P}_{n}$ be the set of all $p$-multipartitions of
$n$.

\addtocounter{dfn}{4} \begin{dfn} {\rm (\cite{AM})} Suppose $n\geq
0$. The set $\mathcal{K}_n$ of Kleshchev $p$ multipartitions with
respect to $(q,Q_1,\cdots,Q_p)$ is defined inductively as follows:
\smallskip

(1)
$\mathcal{K}_0:=\Bigl\{\underline{\emptyset}:=\bigl(\underbrace{\emptyset,\cdots,
\emptyset}_{p}\bigl)\Bigr\}$;\smallskip

(2)
$\mathcal{K}_{n+1}:=\Bigl\{\mu\in\mathcal{P}_{n+1}\Bigm|\text{$\lam\overset{x}
{\twoheadrightarrow}\mu$ for some $\lam\in\mathcal{K}_n$ and some
$x$}\Bigr\}$.
\end{dfn}

The {\it Kleshchev's good lattice} with respect to
$(q,Q_1,\cdots,Q_p)$ is, by definition, the infinite graph whose
vertices are the Kleshchev $p$-multipartitions with respect to
$(q,Q_1,\cdots,Q_p)$ and whose arrows are given by $$
\text{$\lam\overset{x}{\twoheadrightarrow}\mu$\quad$\Longleftrightarrow$\quad
$\lam$ is obtained from $\mu$ by removing a good $x$-node}. $$

Now we can state Ariki's remarkable result (\cite{A2}), that is,

\addtocounter{thm}{3} \begin{thm} {\rm (\cite{A2})} \label{Ariki}
With the above notations, we have that $\td_{\bbQ}^{\lam}\neq 0$ if
and only if $\lam$ is a Kleshchev $p$-multipartition with respect to
$(q,Q_1,\cdots,Q_p)$.
\end{thm}

\medskip\bigskip

\section{Classification of simple $\H_q(p,p,n)$-modules}

In this section, we shall first review some known results on the
classification of simple modules over $\H_q(p,p,n)$. Then we shall
state the first two main results in this paper, which give a
classification of simple modules over $\H_q(p,p,n)$ in the
non-separated cases. The proof will be given in Section 4 and
Section 5.\smallskip

Recall that $\bQ=\bigl\{1,\eps,\cdots,\eps^{p-1}\bigr\}$, where
$\eps$ is a fixed primitive $p$-th root of unity in $K$. Let
$\bbQ$ be an ordered $p$-tuple which is obtained by fixing an
order on $\bQ$. Let ${K}_n=\bigl\{\lam\in\mathcal{P}_{n}\bigm|
\td_{\bbQ}^{\lam}\neq 0\bigr\}$. The automorphism $\sig$
determines uniquely an automorphism $\HH$ of ${K}_n$ such that
$\bigl(\td_{\bbQ}^{\lam}\bigr)^{\sig}\cong\td_{\bbQ}^{\HH(\lam)}$.
Clearly, $\HH^p=\id$. In particular, we get an action of the
cyclic group $C_p$ on ${K}_n$ given as follows:
$\td_{\bbQ}^{\sig^{k}\cdot\lam}\cong
(\td_{\bbQ}^{\lam})^{\sig^k}$. Let $\sim_{\sig}$ be the
corresponding equivalence relation on ${K}_n$. That is,
$\lam\sim_{\sig}\mu$ if and only if $\lam=g\cdot\mu$ for some
$g\in C_p$. For each $\lam\in{K}_n/{\sim_{\sig}}$, let
${C}_{\lam}$ be the stabilizer of $\lam$ in $C_p$. Then
${C}_{\lam}$ is a cyclic subgroup of $C_p$ with order
$|{C}_{\lam}|$. Clearly $|{C}_{\lam}|\mid p$. We define
$$
{K}_n(0):=\bigl\{\lam\in{K}_n/{\sim_{\sig}}\bigm|
{C}_{\lam}=1\bigr\},\quad
{K}_n(1):=\bigl\{\lam\in{K}_n/{\sim_{\sig}}\bigm| {C}_{\lam}\neq 1
\bigr\}. $$

\begin{lem} \label{thm31}{\rm (\cite[(5.4),(5.5),(5.6)]{Hu3})} Suppose $\H_{q}(p,p,n)$ is split over $K$.
\smallskip

1) Let $\td_{\bbQ}^{\lam}$ be any given irreducible
$\H_{q}(p,n)$-module and $D$ be an irreducible
$\H_{q}(p,p,n)$-submodule of $\td_{\bbQ}^{\lam}$. Let $d$ be the
smallest positive integer such that $D\cong DT_{0}^{d}$. Suppose
$1\leq d<p$. Then $k:=p/d$ is the smallest positive integer such
that $\td_{\bbQ}^{\lam}\cong(\td_{\bbQ}^{\lam})^{\sig^k}$, and $$
\td_{\bbQ}^{\lam}\downarrow_{\H_{q}(p,p,n)}\cong D\oplus
DT_0\oplus\cdots\oplus DT_0^{d-1}. $$

2) The set $
\Bigl\{\td_{\bbQ}^{\lam}\downarrow_{\H_{q}(p,p,n)}\Bigm|\lam\in{K}_n(0)\Bigr\}
\bigcup
\Bigl\{D^{\lam,0},D^{\lam,1},\cdots,D^{\lam,|{C}_{\lam}|-1}\Bigm|
\lam\in{K}_n(1)\Bigr\} $ forms a complete set of pairwise
non-isomorphic simple $\H_{q}(p,p,n)$-modules, where for each
$\lam\in{K}_n(1)$, $D^{\lam,0}$ is an irreducible $\H_q(p,p,n)$
submodule of $\td_{\bbQ}^{\lam}$, and
$D^{\lam,i}=\bigl(D^{\lam,0}\bigr)^{\tau^i}$ for
$i=0,1,\cdots,|{C}_{\lam}|-1$.
\end{lem}

Therefore, the problem on classifying simple
$\H_{q}(p,p,n)$-modules reduces to the problem on determining the
automorphism $\HH$.

\addtocounter{dfn}{1}
\begin{dfn} {\rm (\cite{A5})}\label{sep} We refer to the condition $f_{p,n}(q,\eps)\in K^{\times}$ as the separation condition. Say that
we are in the separated case if the separation condition is
satisfied, otherwise we are in the non-separated case.
\end{dfn}

In the separated case, the classification of all the simple
$\H_q(p,p,n)$-modules is known by the results in \cite{Hu3}. In
particular, we have that

\addtocounter{lem}{1} \begin{lem} {\rm (\cite{Hu3})} Suppose
$f_{p,n}(q,\eps)\neq 0$ in $K$. Then $\H_{q}(p,p,n)$ is split over
$K$, and for any $\lam=(\lam^{(1)},\cdots,\lam^{(p)})\in{K}_n$,
$$ \HH(\lam)=\lam[1]:=(\lam^{(2)},\lam^{(3)},\cdots,\lam^{(p)},\lam^{(1)}).
$$\end{lem}

The above result generalizes the corresponding results in
\cite[(3.6),(3.7)]{P} and \cite{Hu1} on the Iwahori--Hecke algebra
of type $D_n$. So it remains to deal with the case when
$f_{p,n}(q,\eps)=0$, i.e., the separated case. In this case, by
Theorem \ref{Ariki}, $K_n=\mathcal{K}_n$. Henceforth we identify
$K_n$ with $\mathcal{K}_n$ without further comments. Note that in
the special case where $p=2$, $q$ must be a primitive $(2\l)$-th
root of unity for some positive integer $\l$ and the classification
of all the simple modules is also known by the results in
\cite{Hu2}. In that case, $\HH$ is an involution. We have that

\begin{lem} {\rm (\cite[(1.4)]{Hu2}, \cite[Appendix]{Hu4})} \label{H2} Suppose $q$ is a primitive
$2\l$-th root of unity for some positive integer $\l$. Suppose
$\ch K\neq 2$ and $\H_q(D_n)$ is split over $K$. Let
$\lam\in\mathcal{K}_n$ be a Kleshchev bipartition (i.e.,
$2$-multipartition) of $n$ with respect to $(\sqrt[2\l]{1},1,-1)$,
and let
$$
\underline{\emptyset}\overset{r_1}{\twoheadrightarrow}\cdot
\overset{r_2}{\twoheadrightarrow}\cdot \cdots\cdots
\overset{r_n}{\twoheadrightarrow}\lam $$ be a path from
$\underline{\emptyset}$ to $\lam$ in Kleshchev's good lattice with
respect to $(\sqrt[2\l]{1},1,-1)$. Then, the sequence $$
\underline{\emptyset}\overset{\l+r_1}{\twoheadrightarrow}\cdot
\overset{\l+r_2}{\twoheadrightarrow}\cdot \cdots\cdots
\overset{\l+r_n}{\twoheadrightarrow}\cdot
$$
also defines a path in Kleshchev's good lattice with respect to
$(\sqrt[2\l]{1},1,-1)$, and it connects $\underline{\emptyset}$ to
$\HH(\lam)$.
\end{lem}
Note that the above description of the involution $\HH$ bears much
resemblance with Kleshchev's description of the well-known
Mullineux involution (see \cite{Kl}, \cite{B}).\smallskip

The first one of the two main results in this section is a direct
generalization of Lemma \ref{H2} to the case of the cyclotomic Hecke
algebra $\H_q(p,p,n)$. By assumption, $K$ is a field which contains
a primitive $p$-th root of unity $\eps$. In particular, $\ch K$ is
coprime to $p$. Now $f_{p,n}(q,\eps)=0$ implies that
$\langle\eps\rangle\cap\langle q\rangle\neq\{1\}$, where
$\langle\eps\rangle$ (resp., $\langle q\rangle$) denotes the
multiplicative subgroup generated by $\eps$ (resp., by $q$). Let
$1\leq k< p$ be the smallest positive integer such that
$\eps^{k}\in\langle q\rangle$. In particular, $k|p$. Suppose
$\eps^{k}=q^{\ell}$, where $\ell>0$. Let $d=p/k$. then
$q^{\ell}=\eps^{k}$ is a primitive $d$-th root of unity. Hence $q$
is a primitive $(d\ell_1)$-th root of unity for some positive
integer $\ell_1$ with $\ell_1|\ell$. We need the following result
from number theory.

\begin{lem} \label{number} Let $K$ be a field which contains
a primitive $p$-th root of unity $\eps$. Suppose $p=dk$, where
$p,d,k\in\mathbb{N}$. $\xi\in K$ is a primitive $d$-th root of
unity. Then there exists a primitive $p$-th root of unity
$\zeta\in K$ such that $\zeta^k=\xi$.
\end{lem}

\noindent {Proof:} \,Clearly $\eps^k$ is a primitive $d$-th root
of unity in $K$. Indeed, the set
$\bigl\{1,\eps,\eps^2,\cdots,\eps^{p-1}\bigr\}$ is the set of all
$p$-th root of unity in $K$, and the set
$\bigl\{1,\eps^k,\eps^{2k},\cdots,\eps^{(d-1)k}\bigr\}$ is the set
of all $d$-th root of unity in $K$. It follows that there exists
some integer $1\leq a<d$ with $(a,d)=1$ and such that
$\xi=\eps^{ak}$.

We write $k=k'k''$, where $(k',d)=1$ and any prime factor of $k''$
is a factor of $d$. By the Chinese Remainder Theorem,
$\Z_{k'd}\cong \Z_{k'}\times\Z_{d}$. Then we can find $j$ such
that $a+jd\equiv 1\pmod{k'}$. In particular, $(a+jd,k')=1$. It
follows that $(a+jd,k)=1$, and hence $(a+jd,p)=1$ (because
$(a,d)=1$). Now $\zeta:=\eps^{a+jd}$ is a primitive $p$-th root of
unity, and $\zeta^k=\eps^{ak+jdk}=\eps^{ak+jp}=\xi$. This
completes the proof of the lemma. \hfill\qed\medskip

We return to our discussion above Lemma \ref{number}. Note that $q$
is a primitive $(d\ell_1)$-th root of unity implies that
$q^{\ell_1}$ is a primitive $d$-th root of unity. By Lemma
\ref{number}, we can always find a primitive $p$-th root of unity
$\widetilde{\eps}$ such that $(\widetilde{\eps})^{k}=q^{\ell_1}$.
Since both $\eps$ and $\widetilde{\eps}$ are primitive $p$-th roots
of unity, it follows that there exist integers $i, j$, such that
$\eps=(\widetilde{\eps})^i$ and $\widetilde{\eps}=\eps^j$. In
particular, $\eps^m\in\langle q\rangle$ if and only if
$(\widetilde{\eps})^m\in\langle q\rangle$ for any integer $m\geq 0$.
Therefore, $1\leq k<p$ is also the smallest positive integer such
that $(\widetilde{\eps})^{k}\in\langle q\rangle$. Replacing $\eps$
by $\widetilde{\eps}$ (which makes the Hecke algebra $\H_{q}(p,n)$
itself unchanged) if necessary, we can assume without loss of
generality that $\ell=\ell_1$. {\it Henceforth, we fix such $\eps$.
Therefore, $q$ is a primitive $d\ell$-th root of unity,
$q^{\ell}=\eps^{k}$ is a primitive $d$-th root of unity, and $1\leq
k<p$ is still the smallest positive integer such that
$\eps^{k}\in\langle q\rangle$.}

For integer $i=1,2,\cdots, k$, we set
$\bbQ_i=(\veps^{i-1},\veps^{k+i-1},\cdots, \veps^{(d-1)k+i-1})$.
Then $\bQ=\bQ_1\sqcup\cdots\sqcup\bQ_k$ is a partition of the
parameter set $\bQ$ into different $q$-orbits. Let
$\bbQ=\bigl(\bbQ_1,\bbQ_2,\cdots,\bbQ_{k}\bigr)$ (concatenation of
ordered tuples). For each $p$-multipartition
$\lam=(\lam^{(1)},\cdots,\lam^{(p)})$ of $n$, we write
$$
\lam^{[i]}=(\lam^{((i-1)d+1)},\lam^{((i-1)d+2)},\cdots,\lam^{(id)}),\,\,\,\text{for}\,\,i=1,2,\cdots,k.
$$ and we use $\theta$ to denote the map $\lam\mapsto
(\lam^{[1]},\cdots,\lam^{[k]})$. Now we can state the first main
result, which deals with the case where $k=1$ and $K=\mathbb{C}$.

\addtocounter{thm}{5}
\begin{thm} \label{main1} Suppose that $K=\mathbb{C}$, and $q,\eps\in K$ such that
$\eps=q^{\ell}$ is a primitive $p$-th root of unity and $q$ is a
primitive $(p\ell)$-th root of unity. Recall our definition of $\HH$
in the second paragraph of this section. Let $\lam\in\mathcal{K}_n$
be a Kleshchev $p$-multipartition of $n$ with respect to
$(q,1,\veps,\cdots,\veps^{p-1})$, and let
$$ \underline{\emptyset}\overset{r_1}{\twoheadrightarrow}\cdot
\overset{r_2}{\twoheadrightarrow}\cdot \cdots\cdots
\overset{r_n}{\twoheadrightarrow}\lam $$ be a path from
$\underline{\emptyset}$ to $\lam$ in Kleshchev's good lattice with
respect to $(q,1,\veps,\cdots,$ $\veps^{p-1})$. Then, the sequence
$$ \underline{\emptyset}\overset{\l+r_1}{\twoheadrightarrow}\cdot
\overset{\l+r_2}{\twoheadrightarrow}\cdot \cdots\cdots
\overset{\l+r_n}{\twoheadrightarrow}\cdot
$$
also defines a path in Kleshchev's good lattice with respect to
$(q,1,\veps,\cdots,$ $\veps^{p-1})$, and it connects
$\underline{\emptyset}$ to $\HH(\lam)$.
\end{thm}

The proof of Theorem \ref{main1} will be given in Section 4. Here we
give an example. Suppose that $n=p=3, q=\sqrt[6]{1}$ and
$\veps=q^2$. Then the following are all Kleshchev
$3$-multipartitions with respect to $(q,1,\veps,\veps^2)$ of $3$
$$\begin{aligned}
&\bigl(\emptyset, \emptyset, (1^3)\bigr),\,\,\bigl(\emptyset,
\emptyset, (2,1)\bigr),\,\, \bigl(\emptyset, \emptyset,
(3)\bigr),\,\,\bigl(\emptyset, (1), (1^2)\bigr),\,\,
\bigl(\emptyset, (1), (2)\bigr),\\
 &\bigl(\emptyset, (1^2),
(1)\bigr),\,\,\bigl(\emptyset, (1^3),
\emptyset\bigr),\,\,\bigl(\emptyset,
(2),(1)\bigr),\,\,\bigl(\emptyset, (2,1), \emptyset\bigr),
\,\,\bigl((1), \emptyset, (1^2)\bigr),\\
 &\bigl((1), \emptyset,
(2)\bigr),\,\,\bigl((1), (1), (1)\bigr), \,\,\bigl((1), (1^2),
\emptyset),\,\,\bigl((1), (2), \emptyset\bigr),\,\, \bigl((1^2),
\emptyset, (1)\bigr),\\
 &\bigl((1^2),
(1),\emptyset\bigr),\,\,\bigl((2), \emptyset, (1)\bigr),
\,\,\bigl((2), (1), \emptyset\bigr),\,\,\bigl((2,1), \emptyset,
\emptyset),\end{aligned}
$$
and the automorphism $\HH$ is given by $$
\begin{matrix}
\bigl(\emptyset, \emptyset, (1^3)\bigr)&\longmapsto & \bigl((1^2),
\emptyset, (1)\bigr)&\longmapsto &\bigl(\emptyset, (1^3),
\emptyset\bigr),\\
\bigl(\emptyset, \emptyset, (2,1)\bigr)&\longmapsto &\bigl((2,1),
\emptyset, \emptyset\bigr)&\longmapsto &\bigl(\emptyset, (2,1),
\emptyset\bigr),\\
\bigl(\emptyset, \emptyset, (3)\bigr)&\longmapsto &\bigl((2), (1),
\emptyset\bigr)&\longmapsto &\bigl(\emptyset, (2),
(1)\bigr),\\
\bigl(\emptyset, (1), (1^2)\bigr)&\longmapsto & \bigl((1),
\emptyset, (2)\bigr)&\longmapsto &\bigl((1), (1^2),
\emptyset\bigr),\\
\bigl(\emptyset, (1), (2)\bigr)&\longmapsto &\bigl((2), \emptyset,
(1)\bigr)&\longmapsto &\bigl((1), (2),
\emptyset\bigr),\\
\bigl(\emptyset, (1^2), (1)\bigr)&\longmapsto &\bigl((1),
\emptyset, (1^2)\bigr)&\longmapsto &\bigl((1^2), (1),
\emptyset\bigr),\\
\bigl((1), (1), (1)\bigr)&\longmapsto &\bigl((1), (1), (1)\bigr)&
\longmapsto &\bigl((1), (1), (1)\bigr).\\
\end{matrix}
$$

In view of Theorem \ref{main1} and the discussion above it, it
remains to consider the case where $k>1$. Now we suppose $k>1$.
Recall our assumption, that is, $q$ is a primitive $(d\ell)$-th root
of unity, $q^{\ell}=\eps^{k}$ is a primitive $d$-th root of unity,
and $1\leq k<p$ is the smallest positive integer such that
$\eps^{k}\in\langle q\rangle$. Let $\eps'=\eps^k$, which is a
primitive $d$-th root of unity in $\mathbb{C}$. For any integer
$n'\geq 1$, let $\H_{q}(d,n')$ be the Ariki--Koike algebra (over
$\mathbb{C}$) with parameters
$\bigl\{q,1,\eps',(\eps')^2,\cdots,(\eps')^{d-1}\bigr\}$ and size
$n'$. Let $\sig'$ be the nontrivial $\mathbb{C}$-algebra
automorphism of $\H_{q}(d,n')$ which is defined on generators by
$\sig(T_0)=\eps' T_0, \sig(T_i)=T_i$, for $i=1,2,\cdots, n'-1$. Let
$\mathcal{K}'_{n'}$ be the set of Kleshchev $d$-multipartitions of
$n'$ with respect to $(q,1,\eps',(\eps')^2,\cdots,(\eps')^{d-1})$.
Let $\bbQ':=(1,\eps',(\eps')^2,\cdots,(\eps')^{d-1})$. Note that all
the parameters in $\bbQ'$ are in a single $q$-orbit. For any
$d$-multipartition $\lam'$ of $n'$, by Theorem \ref{Ariki},
$\td_{\bbQ'}^{\lam'}\neq 0$ if and only if
$\lam'\in\mathcal{K}'_{n'}$. The automorphism $\sig'$ determines
uniquely an automorphism $\HH'$ of $\mathcal{K}'_{n'}$ such that
$\bigl(\td_{\bbQ'}^{\lam'}\bigr)^{\sig'}\cong\td_{\bbQ'}^{\HH'(\lam')}$.
Clearly, $(\HH')^d=\id$. Note that $\eps'=\eps^k=q^{\ell}$. Hence we
are in a position to apply Theorem \ref{main1} with $p$ replaced by
$d$, $n$ replaced by $n'$ and $\eps$ replaced by $\eps'$. We get
that

\addtocounter{cor}{6}
\begin{cor} \label{maincor} Let $\lam'\in\mathcal{K}'_{n'}$ be a Kleshchev $d$-multipartition of
$n'$ with respect to $(q,1,\veps',\cdots,{(\veps')}^{d-1})$, and let
$$ \bigl(\underbrace{\emptyset,\cdots,
\emptyset}_{d}\bigl)=\underline{\emptyset}\overset{r_1}{\twoheadrightarrow}\cdot
\overset{r_2}{\twoheadrightarrow}\cdot \cdots\cdots
\overset{r_{n'}}{\twoheadrightarrow}\lam' $$ be a path from
$\underline{\emptyset}$ to $\lam'$ in Kleshchev's good lattice with
respect to $(q,1,\veps',\cdots,$ ${(\veps')}^{d-1})$. Then, the
sequence
$$ \bigl(\underbrace{\emptyset,\cdots,
\emptyset}_{d}\bigl)=\underline{\emptyset}\overset{\l+r_1}{\twoheadrightarrow}\cdot
\overset{\l+r_2}{\twoheadrightarrow}\cdot \cdots\cdots
\overset{\l+r_{n'}}{\twoheadrightarrow}\cdot
$$
also defines a path in Kleshchev's good lattice with respect to
$(q,1,\veps',\cdots,$ ${(\veps')}^{d-1})$, and it connects
$\underline{\emptyset}$ to $\HH'(\lam')$.
\end{cor}

Now we can state the second main result, which deals with the case
where $k>1$ and $K=\mathbb{C}$.

\addtocounter{thm}{1}
\begin{thm}\label{main2} Suppose that $K=\mathbb{C}$, $q,\eps\in K$ such
that $\eps$ is a primitive $p$-th root of unity, $\eps^k=q^{\ell}$
is a primitive $d$-th root of unity and $q$ is a primitive
$(d\ell)$-th root of unity, and $1\leq k<p$ is the smallest positive
integer such that $\eps^{k}\in\langle q\rangle$. Let
$\lam\in\mathcal{K}_n$ be a Kleshchev $p$-multipartition of $n$ with
respect to $(q,\bbQ)$, where
$\bbQ=\bigl(\bbQ_1,\bbQ_2,\cdots,\bbQ_{k}\bigr)$ (concatenation of
ordered tuples). Then
$$\theta\Bigl(\HH(\lam)\Bigr)=\Bigl(\HH'(\lam^{[k]}),\lam^{[1]},\cdots,\lam^{[k-1]}\Bigr),
$$
where $|\lam^{[k]}|=n'$, $\HH'$ is as defined in Corollary
\ref{maincor} and the righthand side of the above equality is
understood as concatenation of ordered tuples.
\end{thm}

The proof of Theorem \ref{main2} will be given in Section 5.

\begin{thm} \label{main3} Both Theorem \ref{main1} and Theorem \ref{main2}
remain true if we replace $\mathbb{C}$ by any field $K$ such that
$\H_q(p,p,n)$ is split over $K$ and $K$ contains primitive $p$-th
root of unity.
\end{thm}

\noindent {Proof:} \,In the case where $p=2$, this is proved in the
appendix of \cite{Hu4} (which is essentially an argument due to S.
Ariki). In general, this can still be proved by using the same
argument as in the appendix of \cite{Hu4}.\hfill\qed
\medskip

Note that we have proved in \cite[Theorem 5.7]{Hu3} that in the
separated case $\H_q(p,p,n)$ is always split over $K$ whenever $K$
contains primitive $p$-th root of unity. It would be interesting to
know if this is still true in the non-separated case.
\medskip\bigskip

\noindent {\bf Remark 3.10}\,\,\,Let $\bbQ:=(Q_1,\cdots,Q_p)$ be an
arbitrary permutation of $(1,\eps,$ $\cdots,\eps^{p-1})$. We
redefine the residue of a node $\gamma=(a,b,c)$ to be
$\rres(\gamma):=q^{b-a}Q_c$. The notions of good nodes and Kleshchev
multipartitions are defined in a similar way as before. For any
$r\in K$, we use $\lam\overset{r}{\twoheadrightarrow}\mu$ to
indicate $\lam$ is obtained from $\mu$ by removing a good $r$-node.
One can easily deduce from Theorem \ref{main1} and Theorem
\ref{main2} the following description of $\HH$: let
$\lam\in\mathcal{K}_n$ be a Kleshchev $p$-multipartition of $n$ with
respect to $(q,Q_1,\cdots,Q_p)$, and let
$\underline{\emptyset}\overset{r_1}{\twoheadrightarrow}\cdot
\overset{r_2}{\twoheadrightarrow}\cdot \cdots\cdots
\overset{r_n}{\twoheadrightarrow}\lam $ be a path from
$\underline{\emptyset}$ to $\lam$ in Kleshchev's good lattice with
respect to $(q,Q_1,\cdots,$ $Q_p)$. Then, the sequence
\addtocounter{equation}{10}\begin{equation}\underline{\emptyset}\overset{\eps
r_1 }{\twoheadrightarrow}\cdot \overset{\eps r_2
}{\twoheadrightarrow}\cdot \cdots\cdots \overset{\eps r_n
}{\twoheadrightarrow}\cdot\label{311}\end{equation} also defines a
path in Kleshchev's good lattice with respect to
$(q,Q_1,\cdots,Q_p)$, and it connects $\underline{\emptyset}$ to
$\HH(\lam)$.
\medskip

\noindent {\bf Remark 3.12}\,\,\,In the case of $G(r,p,n)$, one can
still define the automorphisms $\sig, \tau$ and $\HH$, and the
problem of classifying simple $\H_q(r,p,n)$-modules can still be
reduced to the determination of $\HH$. Note that although we deal
with the $G(p,p,n)$ case in this paper, it is not difficult to
generalize Theorem \ref{main1} and Theorem \ref{main2} in this paper
to the $G(r,p,n)$ case by using the same argument as well as the
Morita equivalence theorem (\cite[Theorem 1.1]{DM}) of
Dipper--Mathas for Ariki--Koike algebras\footnote{We do not deal
with the general $G(r,p,n)$ case here because this paper has already
been cited in Ariki's book \cite[{[}cyclohecke12{]}]{A3} on the one
hand; and on the other hand, the proof in the $G(r,p,n)$ case needs
nothing more than the proof in the $G(p,p,n)$ case except some
sophisticated notations.}. In the preprint \cite{GJ}, Genet and
Jacon give a parameterization of simple $\H_q(r,p,n)$-modules when
$q$ is a root of unity by using combinatorics of FLOTW partitions.
Our $\sig,\tau, \HH$ are denoted by $f,g,\tau$ in their paper. They
give a characterization of $\tau$ in terms of $\omega$ and a
bijection $\kappa$ between the set of Kleshchev multipartitions and
the set of FLOTW partitions (see \cite[Proposition 2.9, Line3--4 in
Page 14]{GJ}). Note that our description of $\HH$ in (\ref{311}) is
actually a statement about the crystal graph, and although \cite{GJ}
uses JMMO's Fock space (\cite{JMM}, \cite{FLO}) and we use Hayashi's
Fock space (\cite{Ha}, \cite{AM}), the two crystals provided by
lattice of Kleshchev multipartitions and by lattice of FLOTW
partitions respectively are isomorphic to each other. It follows
that the description of $\HH$ in (\ref{311}) in the context of
Kleshchev's good lattices should also be valid in the context of
FLOTW's good lattices. Therefore, there is a second approach which
can be used to generalize Theorem \ref{main1} and Theorem
\ref{main2} in this paper to the $G(r,p,n)$ case. The details are
given in \cite{Hu5}. The main idea is to derive from the setting of
FLOTW partitions and \cite[Proposition 2.9, Line3--4 in Page
14]{GJ}) a description of $\HH$ like (\ref{311}) (see
\cite[(4.3.1)]{Ja0} for the special case when $p=r$ and
$\eps=q^{\l}$). In the case where $p=r$, the parameterization of
simple modules obtained in \cite{GJ} is the same as \cite[Theorem
5.7]{Hu3} in the separated case. It is worthwhile and interesting to
establish a direct connection between the two parameterizations of
irreducible representations.
\medskip

\section{Proof of Theorem \ref{main1}}

In this section, we shall give the proof of Theorem \ref{main1}. It
turns out that the proof of Theorem \ref{main1} is a direct
generalization of the proof of \cite[(1.5)]{Hu2}. Throughout this
section, we keep the same assumptions and notations as in Theorem
\ref{main1}. That is, $K=\mathbb{C}$,\, $q,\eps\in \mathbb{C}$ be
such that $\eps=q^{\ell}$ is a primitive $p$-th root of unity and
$q$ is a primitive $(p\ell)$-th root of unity.\smallskip

Let $v$ be an indeterminate over $\mathbb{Q}$. Let $\mathfrak{h}$
be a $(p\ell+1)$-dimensional vector space over $\mathbb{Q}$ with
basis $\bigl\{h_0,h_1,\cdots,h_{p\ell-1},d\bigr\}$\footnote{The
readers should not confuse the element $d$ here with the integer
$d$ we used before.}. Denote by $\bigl\{
\Lambda_0,\Lambda_1,\cdots,$ $\Lambda_{p\ell-1},\delta\bigr\}$ the
corresponding dual basis of $\mathfrak{h}^{\ast}$, and we set
$\alpha_i=2\Lambda_{i}-\Lambda_{i-1}-\Lambda_{i+1}+\delta_{i,0}\delta$
for $i\in\Z/p\ell\Z$. The weight lattice is
$P=\Z\Lambda_0\oplus\cdots\oplus\Z\Lambda_{p\ell-1}\oplus
\Z\delta$, its dual is $P^{\vee}=\Z h_0\oplus\cdots\oplus\Z
h_{p\ell-1}\oplus\Z d$. Assume that the $p\ell\times p\ell$ matrix
$\bigl(\langle\alpha_{i}, h_{j}\rangle\bigr)$ is just the
generalized Cartan matrix associated to $\ksl_{p\ell}$. The
quantum affine algebra $U_v({\ksl}_{p\ell})$ is by definition the
$\mathbb{Q}(v)$-algebra with $1$ generated by elements $E_i, F_i,
\,\,i\in\Z/p\ell\Z$\,\, and $K_{h},\,\,h\in P^{\vee}$, subject to
the relations $$\begin{matrix}
K_{h}K_{h'}=K_{h+h'}=K_{h'}K_{h},\,\,\,K_0=1,\\[4pt]
K_{h}E_{j}=v^{\langle\alpha_{j},h\rangle}E_{j}K_{h},\,\,\,
K_{h}F_{j}=v^{-\langle\alpha_{j},h\rangle}F_{j}K_{h},\\[4pt]
E_{i}F_{j}-F_{j}E_{i}=\delta_{ij}\frac{K_{h_i}-K_{-h_i}}{v-v^{-1}},\\[4pt]
\sum\limits_{k=0}^{1-\langle\alpha_{i},h_j\rangle}(-1)^{k}\begin{bmatrix}
1-\langle\alpha_{i},h_j\rangle\\[2pt]
k\end{bmatrix} E_{i}^{1-\langle\alpha_{i},h_j\rangle-k}
E_{j}E_{i}^{k}=0,\,\,\,(i\neq j),\\[4pt]
\sum\limits_{k=0}^{1-\langle\alpha_{i},h_j\rangle}(-1)^{k}
\begin{bmatrix}
1-\langle\alpha_{i},h_j\rangle\\[2pt]
k\end{bmatrix} F_{i}^{1-\langle\alpha_{i},h_j\rangle-k}
F_{j}F_{i}^{k}=0,\,\,\,(i\neq j),
\end{matrix}
$$
where $$ [k]:=\frac{v^{k}-v^{-k}}{v-v^{-1}},\quad
[k]!:=[k][k-1]\cdots[1],\quad
\begin{bmatrix}
m\\[2pt]
k\end{bmatrix}:=\frac{[m]!}{[m-k]![k]!}.
$$
It is a Hopf algebra with comultiplication given by $$
\Delta(K_h)=K_{h}\otimes K_{h},\,\,\, \Delta(E_i)=E_{i}\otimes
1+K_{-h_i}\otimes E_{i},\,\,\, \Delta(F_{i})=F_{i}\otimes
K_{h_i}+1\otimes F_{i}. $$

Let $\mathcal{P}^{(1)}$ be the set of all partitions. Let
$\mathcal{P}:=\sqcup_{n\geq 0}\mathcal{P}_n$, the set of all
$p$-multipartitions. For each integer $j$ with $0\leq j<p\ell$,
there is a {\it level $1$ Fock space}
$\mathcal{F}^{(1)}(\Lambda_j) $, which is defined as follows. As a
vector space,  $$
\mathcal{F}^{(1)}(\Lambda_j):=\bigoplus_{\lam\in\mathcal{P}^{(1)}}\mathbb{Q}(v)\lam,
$$ and the algebra $U_v({\ksl}_{p\ell})$ acts on
$\mathcal{F}^{(1)}(\Lambda_j)$ by
$$\begin{matrix}
K_{h_i}\lam=v^{N_i(\lam)}\lam,&K_{d}\lam=v^{-N_d(\lam)}\lam,\\
E_{i}\lam=\sum_{\nu\overset{i}{\rightarrow}\lam}v^{-N_i^{r}(\nu,\lam)}\nu,&
F_{i}\lam=\sum_{\lam\overset{i}{\rightarrow}\mu}v^{N_i^{l}(\lam,\mu)}\mu,
\end{matrix}
$$
where $i\in\Z/p\ell\Z$, $\lam\in\mathcal{P}^{(1)}$, and
$$\begin{aligned}
N_{i}(\lam)&=\#\Bigl\{\mu\bigm|\lam\overset{i}{\rightarrow}\mu\Bigr\}-
\#\Bigl\{\nu\bigm|\nu\overset{i}{\rightarrow}\lam\Bigr\},\\
N_{i}^{r}(\nu,\lam)&=\sum_{\gamma\in
[\lambda]\setminus[\nu]}\biggl(\#\Bigl\{\gamma'\Bigm|
\begin{matrix}
\text{$\gamma'$ an addable}\\
\text{$i$-node of $\lam$ above $\gamma$}\end{matrix}\Bigr\}-\\
&\qquad\qquad\quad\#\Bigl\{\gamma'\Bigm|\begin{matrix}
\text{$\gamma'$ a removable}\\
\text{$i$-node of $\nu$ above $\gamma$}\end{matrix}\Bigr\}\biggr),\\
N_{i}^{l}(\lam,\mu)&=\sum_{\gamma\in
[\mu]\setminus[\lam]}\biggl(\#\Bigl\{\gamma'\Bigm|
\begin{matrix}
\text{$\gamma'$ an addable}\\
\text{$i$-node of $\mu$ below $\gamma$}\end{matrix}\Bigr\}-\\
&\qquad\qquad\quad\#\Bigl\{\gamma'\Bigm|\begin{matrix}
\text{$\gamma'$ a removable}\\
\text{$i$-node of $\lam$ below
$\gamma$}\end{matrix}\Bigr\}\biggr),\end{aligned}
$$
and
$N_{d}(\lam):=\#\bigl\{\gamma\in[\lam]\bigm|\res(\gamma)=0\bigr\}$,
and here we should use the following definition of residue,
namely, the node in the $a$th row and the $b$th column of $\lam$
is filled out with the residue $b-a+j\in{\Z/p\ell\Z}$.\smallskip

Let
$\Lambda:=\Lambda_0+\Lambda_{\l}+\Lambda_{2\l}+\cdots+\Lambda_{(p-1)\l}$.
Replacing $\mathcal{P}^{(1)}$ by $\mathcal{P}$, and using the
definition of residue given in (\ref{res}) for multipartition, we
get a {\it level $p$ Fock space}
$$ \mathcal{F}(\Lambda):=\bigoplus_{\lam\in\mathcal{P}}\mathbb{Q}(v)\lam. $$
As a vector space, we have
$\mathcal{F}(\Lambda)\cong\mathcal{F}^{(1)}(\Lambda_0)\otimes\cdots\otimes
\mathcal{F}^{(1)}(\Lambda_{(p-1)\l})$,
$\lam\mapsto\lam^{(1)}\otimes\cdots\otimes\lam^{(p)}$. By
\cite[(2.5)]{AM}, it is indeed an $U_v({\ksl}_{p\l})$-module
isomorphism, where the action on the right-hand side is defined
via $\Delta^{(p-1)}:=\underbrace{(\Delta\otimes
1\otimes\cdots\otimes 1)}_{p-1} \cdots\underbrace{(\Delta\otimes
1)}_{2}\underbrace{\Delta}_{1}$.\smallskip

For each $\lam=(\lam^{(1)},\cdots,\lam^{(p)})\in\mathcal{P}$, we
define
$\widehat{\lam}=(\lam^{(p)},\lam^{(1)},\cdots,\lam^{(p-1)})$. Let
${\Theta}$ be the automorphism of
$\Bigl(U_v({\ksl}_{p\l})\Bigr)^{\otimes p}$ which is defined by
$\Theta(x_1\otimes\cdots\otimes x_p)=x_2\otimes \cdots\otimes
x_{p}\otimes x_1$ for any $x_1,\cdots,x_p\in U_v({\ksl}_{p\l})$.
We now introduce a different version of level $p$ Fock space
$\widehat{\mathcal{F}}(\Lambda)$. As a vector space,
$\widehat{\mathcal{F}}(\Lambda)=\mathcal{F}(\Lambda)\cong\mathcal{F}^{(1)}(\Lambda_0)\otimes\cdots\otimes
\mathcal{F}^{(1)}(\Lambda_{(p-1)\l})$, while the action (denoted
by ``$\circ$'') of $U_v({\ksl}_{p\l})$ is defined by $$
x\circ(\lam^{(1)}\otimes\cdots\otimes\lam^{(p)}):=\Bigl\{\Theta\Bigl(
\Delta^{(p-1)}(x)\Bigr)\Bigr\}(\lam^{(1)}\otimes\cdots\otimes\lam^{(p)}).
$$

\begin{lem} \label{lm41} The above action defines an integrable representation of the algebra $U_v({\ksl}_{p\l})$ on
$\widehat{\mathcal{F}}(\Lambda)$, such that $$\begin{matrix}
K_{h_i}\circ\lam=v^{N_i(\lam)}\lam,&K_{d}\circ\lam=v^{-N_d(\lam)}\lam,\\
E_{i}\circ\lam=\sum_{\nu\overset{i}{\rightarrow}\lam}
v^{-N_{i+\ell}^{r}(\widehat{\nu},\widehat{\lam})}\nu,&
F_{i}\circ\lam=\sum_{\lam\overset{i}{\rightarrow}\mu}v^{N_{i+\ell}^{l}
(\widehat{\lam},\widehat{\mu})}\mu.\\
\end{matrix} $$
Moreover, the empty $p$-multipartition
$\underline{\emptyset}:=(\underbrace{\emptyset,\cdots,\emptyset}_{p})$
is a highest weight vector of weight
$\Lambda=\Lambda_{0}+\Lambda_{\l}+\cdots+\Lambda_{(p-1)\l}$, and
$\widehat{L}(\Lambda):=U_v({\ksl}_{p\l})\circ\underline{\emptyset}$
is isomorphic to the irreducible highest weight module with
highest weight $\Lambda$.
\end{lem}

\noindent {Proof:} \, For any $x,y\in U_v({\ksl}_{p\l})$ and
$\lam\in\mathcal{P}$, we have that $$\begin{aligned} x\circ
(y\circ \lam)&=\Bigl(\Theta\bigl(\Delta^{p-1}(x)\bigr)\Bigr)
\biggl(\Bigl(\Theta\bigl(\Delta^{p-1}(y)\bigr)\Bigr)\lam\biggr)\\
&=\Bigl(\Theta\bigl(\Delta^{p-1}(x)\bigr)\Theta\bigl(\Delta^{p-1}(y)\bigr)\Bigr)\lam
=\Bigl(\Theta\bigl(\Delta^{p-1}(x)\Delta^{p-1}(y)\bigr)\Bigr)\lam\\
&=\Bigl(\Theta\bigl(\Delta^{p-1}(xy)\bigr)\Bigr)\lam
=(xy)\circ\lam,\end{aligned}
$$
and $1\circ\lam=K_0\circ\lam=\lam$. It follows that the action
``$\circ$" does define a representation of the algebra
$U_v({\ksl}_{p\l})$ on $\widehat{\mathcal{F}}(\Lambda)$. The
remaining part of the lemma can be verified easily.\hfill\qed
\medskip

Therefore, both $\mathcal{F}(\Lambda)$ and
$\widehat{\mathcal{F}}(\Lambda)$) are integrable
$U_v({\ksl}_{p\l})$-modules. They are not irreducible. In both
cases the empty multipartition $\underline{\emptyset}$ is a
highest weight vector of weight
$\Lambda=\Lambda_{0}+\Lambda_{\l}+\cdots+\Lambda_{(p-1)\l}$, and
both ${L}(\Lambda):=U_v({\ksl}_{p\l})\underline{\emptyset}$ and
$\widehat{L}(\Lambda)$ are isomorphic to the irreducible highest
weight module with highest weight $\Lambda$. Let
$U'_v({\ksl}_{p\l})$ be the subalgebra of $U_v({\ksl}_{p\l})$
generated by $E_{i},F_{i},K_{h_i}^{\pm 1}, i\in\Z/p\l\Z$. Let $\#$
be the automorphism of $U'_v({\ksl}_{p\l})$ which is defined on
generators by
$$E_{i}^{\#}:=E_{\l+i},\,\,F_{i}^{\#}:=F_{\l+i},\,\,K_{h_{i}}^{\#}:=K_{h_{i+\l}},
\quad \forall\,i\in\Z/p\l\Z. $$

Now we begin to follow the streamline of the proof of
\cite[(1.5)]{Hu2}. We first recall some basic facts about crystal
bases. Let $\A:=\mathbb{Q}[v,v^{-1}]$, let $A$ be the ring of
rational functions in $\mathbb{Q}(v)$ which do not have a pole at
$0$. Write
$$ \mathcal{F}(\Lambda)_{\A}:=\bigoplus_{\lam\in\mathcal{P}}\A\lam,\quad
\mathcal{F}(\Lambda)_{A}:=\bigoplus_{\lam\in\mathcal{P}}A\lam, $$
and we use similar notations for $\widehat{\mathcal{F}}(\Lambda)$.
Let $U_{\A}$ be the Lusztig--Kostant $\A$-form of
$U_v({\ksl}_{p\l})$. Then by \cite[(2.7)]{AM} we know that both
$\mathcal{F}(\Lambda)_{\A}$ and
$\widehat{\mathcal{F}}(\Lambda)_{\A}$ are $U_{\A}$-modules.

Let $u_{\Lambda}:=\underline{\emptyset}$, the highest weight
vector of weight $\Lambda$ in $\mathcal{F}(\Lambda)$. For each
$i\in\{0,1,\cdots,p\l-1\}$, let $\widetilde{E}_{i},
\widetilde{F}_{i}$ be the Kashiwara operators introduced in
\cite{Kas}. Let $L(\Lambda)_{A}$ be the $A$-submodule of
$\mathcal{F}(\Lambda)_{A}$ generated by all
$\widetilde{F}_{i_1}\cdots\widetilde{F}_{i_k}u_{\Lambda}$ for all
$i_1,\cdots,i_k\in\Z/{p\l}Z$. It is a free $A$-module, and is
stable under the action of $\widetilde{E}_{i}$ and
$\widetilde{F}_{i}$. The set
$$
\mathbb{B}(\Lambda):=\bigl\{\widetilde{F}_{i_1}\cdots\widetilde{F}_{i_k}u_{\Lambda}+
vL(\Lambda)_{A}\bigm|i_{1},\cdots,i_{k}\in\Z/p\l\Z\bigr\}\setminus\bigl\{0\bigr\}
$$ is a basis of $L(\Lambda)_{A}/vL(\Lambda)_{A}$. In \cite{Kas},
the pair $\bigl(L(\Lambda)_{A},\mathbb{B}(\Lambda)\bigr)$ is
called the {\it lower crystal basis} at $v=0$ of
$L(\Lambda)$.\smallskip

Following \cite{Kas}, the {\it crystal graph} of $L(\Lambda)$ is
the edge labelled direct graph whose set of vertices is
$\mathbb{B}(\Lambda)$ and whose arrows are given by $$
\text{$b\overset{i}{\twoheadrightarrow}b'$\quad$
\Longleftrightarrow$\quad $\widetilde{F}_{i}b=b'$\,\,\, for some
$i\in\Z/p\l\Z.$}
$$
It is a remarkable fact (\cite{MM}, \cite[(2.11)]{AM}) that the
crystal graph of $L(\Lambda)$ is exactly the same as the
Kleshchev's good lattice if we use the embedding
$L(\Lambda)\subset \mathcal{F}(\Lambda)$. In particular,
$\mathbb{B}(\Lambda)$ can be identified with
$\mathcal{K}:=\sqcup_{n\geq 0}\mathcal{K}_n$. Henceforth, we fix
such an identification.
\smallskip

Let ``$-$'' be the involutive ring automorphism of
$U_v({\ksl}_{p\l})$ which is defined by $$ \begin{matrix}
\overline{v}:=v^{-1},\quad \overline{K_{h}}:=K_{-h},\,\,(h\in P^{\vee}),\\
\overline{E_i}=E_{i},\quad \overline{F_{i}}:=F_{i},\,\,i=0,1,\cdots,p\l-1. \\
\end{matrix}
$$
This gives rise to an involution (still denoted by ``$-$'') of
$L(\Lambda)$. That is, for $x=P\underline{\emptyset}\in
L(\Lambda)$, we set $\overline{x}:=\overline{P}
\underline{\emptyset}$. By \cite{Kas}, there exists a unique
$\A$-basis $\bigl\{G(\mu)\bigm|\mu\in\mathcal{K}\bigr\}$ of
$L(\Lambda)_{\A}$ such that \begin{enumerate}
\item[(G1)] $G(\mu)\equiv\mu\pmod{vL(\Lambda)_{A}}$,
\item[(G2)] $\overline{G(\mu)}=G(\mu)$.
\end{enumerate}
The basis $\bigl\{G(\mu)\bigr\}_{\mu\in\mathcal{K}}$ is called the
{\bf lower global crystal basis} of $L(\Lambda)$. Let
$\lam\in\mathcal{P}_n, \mu\in\mathcal{K}_n$. Let
$d_{\lam,\mu}:=[\ts_{\bbQ}^{\lam}:\td_{\bbQ}^{\mu}]$, and let
$d_{\lam,\mu}(v)\in\A$ be such that
$G(\mu)=\sum_{\lam}d_{\lam,\mu}(v)\lam$. By a well-known result of Ariki \cite{A1},
$d_{\lam,\mu}(1)=d_{\lam,\mu}$, for any $\lam\in\mathcal{P}_n, \mu\in\mathcal{K}_n$.
\medskip

The following six results can be proved by using exactly the same
arguments as in the proof of
\cite[(2.2),(3.1),(3.2),(3.3),(3.4),(3.5)]{Hu2}.

\begin{lem} \label{thm42} Let $\mathcal{F}(\Lambda)^{\#}$ be the
$U'_v({\ksl}_{p\l})$-module which is obtained from
$\mathcal{F}(\Lambda)$ by twisting the action by the automorphism
$\#$. Let $\phi$ be the linear map
$\mathcal{F}(\Lambda)^{\#}\rightarrow\widehat{\mathcal{F}}(\Lambda)$
which is defined by
$\sum_{\lambda}f_{\lambda}(v)\lambda\mapsto\sum_{\lambda}f_{\lambda}(v)\widehat{\lambda}$.
Then $\phi$ is a $U'_v({\ksl}_{p\l})$-module isomorphism.
Moreover, $\phi\bigl(L(\Lambda)\bigr)=\widehat{L}(\Lambda)$.
\end{lem}

\begin{lem} Let $U_v({\ksl}_{p\l})^{-}$ be the subalgebra of
$U_v({\ksl}_{p\l})$ generated by $F_{i}, i\in\Z/p\l\Z$. Then there
is a $\mathbb{Q}(v)$-linear automorphism (denoted by $\#$) on the
irreducible $U_v({\ksl}_{p\l})$-module $L(\Lambda)$ such that $$
\bigl(\underline{\emptyset}\bigr)^{\#}:=\underline{\emptyset},\,\,
\bigl(Px\bigr)^{\#}:=P^{\#}x^{\#},\quad \forall\,\,P\in
U_v({\ksl}_{p\l})^{-}, x\in L(\Lambda). $$
\end{lem}

\begin{lem} \label{lm45} For any $x\in L(\Lambda)$, we have
$\bigl(\widetilde{F}_{i}x\bigr)^{\#} =\widetilde{F}_{\l+i}x^{\#}$.
\end{lem}

\addtocounter{cor}{4}
\begin{cor} Let $\lam\in\mathcal{K}_n$ be a
Kleshchev $p$-multipartition of $n$ with respect to
$(q,1,\eps,\eps^2,\cdots,\eps^{p-1})$, and let $
\underline{\emptyset}\overset{r_1}{\twoheadrightarrow}\cdot
\overset{r_2}{\twoheadrightarrow}\cdot \cdots\cdots
\overset{r_n}{\twoheadrightarrow}\lam $ be a path from
$\underline{\emptyset}$ to $\lam$ in Kleshchev's good lattice with respect to
$(q,1,\eps,\eps^2,\cdots,\eps^{p-1})$.
Then, the sequence $$
\underline{\emptyset}\overset{\l+r_1}{\twoheadrightarrow}\cdot
\overset{\l+r_2}{\twoheadrightarrow}\cdot \cdots\cdots
\overset{\l+r_n}{\twoheadrightarrow}\cdot
$$
also defines a path in Kleshchev's good lattice with respect to
$(q,1,\eps,\eps^2,\cdots,$ $\eps^{p-1})$. We denote the
endpoint by $\lam^{\#}$.
\end{cor}

\addtocounter{lem}{1}
\begin{lem} For each $\mu\in\mathcal{K}$,
$G(\mu)^{\#}=G\bigl(\mu^{\#}\bigr)$.
\end{lem}

\begin{lem} \label{lm48} For any $x\in L(\Lambda)$, we have
$\phi\bigl(\widetilde{F}_{i}x\bigr)
=\widetilde{F}_{\l+i}\circ\phi(x)$.
\end{lem}

Recall that $\widehat{L}(\Lambda)$ is also an irreducible highest
weight $U_v({\ksl}_{p\l})$-module of highest weight $\Lambda$. For
any $i_{1},\cdots,i_{k}\in\Z/p\l\Z$, it is easy to see
$\widetilde{F}_{i_1}\cdots\widetilde{F}_{i_k}u_{\Lambda}\in
vL(\Lambda)_{A}$ if and only if
$\widetilde{F}_{i_1+\l}\circ\cdots\circ\widetilde{F}_{i_k+\l}\circ
u_{\Lambda}\in v\widehat{L}(\Lambda)_{A}$. By Lemma \ref{thm42}
and Lemma \ref{lm48} we know that the lower global crystal basis
of $\widehat{L}(\Lambda)$ is parameterized by
$\widehat{\mathcal{K}}:=\bigl\{\widehat{\mu}\bigm|\mu\in\mathcal{K}\bigr\}$.
We denote them by
$\bigl\{\widehat{G}(\widehat{\mu})\bigm|\mu\in\mathcal{K}\bigr\}$.
For any $\lam\in\mathcal{P}, \mu\in\mathcal{K}$, let
$\widehat{d}_{\lam,\widehat{\mu}}(v)\in\A$ be such that
$\widehat{G}(\widehat{\mu})=\sum_{\lam}\widehat{d}_{\lam,\widehat{\mu}}(v)\lam$.
The following three results can be proved by using exactly the
same arguments as in the proof of \cite[(3.6),(3.7),(3.8)]{Hu2}.

\addtocounter{cor}{2}
\begin{cor}\label{cor49} For any $\lam\in\mathcal{P},
\mu\in\mathcal{K}$, we have
$\phi\bigl(G(\mu)\bigr)=\widehat{G}(\widehat{\mu})$, and
$d_{\lam,\mu}(v)= \widehat{d}_{\widehat{\lam},\widehat{\mu}}(v)$.
\end{cor}

\addtocounter{lem}{1}
\begin{lem} \label{thm410} Let $\varphi\colon
L(\Lambda)\rightarrow\widehat{\mathcal{F}}(\Lambda)$ be the map
which is defined by $\varphi(x):=\phi\bigl(x^{\#}\bigr)$. Then, if
specialized at $v=1$, the $\varphi$ is the restriction of the
${\ksl}_{p\l}$-module isomorphism
$\mathcal{F}(\Lambda)_{\mathbb{Q}}\rightarrow\widehat{\mathcal{F}}(\Lambda)_{\mathbb{Q}}$
given by $\lam\mapsto\lam$.
\end{lem}

\addtocounter{cor}{1}
\begin{cor}\label{cor411} For any $\lam\in\mathcal{P},
\mu\in\mathcal{K}$, we have
$\varphi\bigl(G(\mu^{\#})\bigr)=\widehat{G}\bigl(\widehat{\mu}\bigr)$
and $d_{\lam,\mu^{\#}}(1)=\widehat{d}_{{\lam},\widehat{\mu}}(1)$.
\end{cor}

\bigskip
\noindent {\bf Proof of Theorem \ref{main1}:}\quad  It suffices to
show that $\HH(\mu)=\mu^{\#}$ for any $\mu\in\mathcal{K}_{n}$. It
is well-known that
$\bigl(\ts_{\bbQ}^{\widehat{\lam}}\bigr)_{\mathbb{C}(v)}\cong
\bigl(\ts_{\bbQ}^{\lam}\bigr)_{\mathbb{C}(v)}^{\sigma}$ (see e.g.,
\cite[(3.7),(5.8)]{Hu3}). Hence in the Grothendieck group of the
category of finite-dimensional $\H_{q}(p,n)$-modules,
$[\ts_{\bbQ}^{\widehat{\lam}}]=[\bigl(\ts_{\bbQ}^{\lam}\bigr)^{\sigma}]$.
By Corollary \ref{cor49}, Corollary \ref{cor411} and Ariki's result \cite{A1}, we deduce that
$$\begin{aligned}
{[}\ts_{\bbQ}^{\widehat{\lam}}:{\td}_{\bbQ}^{\HH(\mu)}{]}
&=[\bigl(\ts_{\bbQ}^{\lam}\bigr)^{\sigma}:\bigl(\td_{\bbQ}^{\mu}
\bigr)^{\sigma}]=[\ts_{\bbQ}^{\lam}:\td_{\bbQ}^{\mu}]=d_{\lam,\mu}\\
&=d_{\lam,\mu}(1)=\widehat{d}_{\widehat{\lam},\widehat{\mu}}(1)
=d_{\widehat{\lam},\mu^{\#}}(1)=d_{\widehat{\lam},\mu^{\#}}
=[\ts_{\bbQ}^{\widehat{\lam}}:\td_{\bbQ}^{\mu^\#}],\end{aligned}
$$ for any $\lam\in\mathcal{P}_n$. Taking $\widehat\lam$ to be
$\HH(\mu)$ or $\mu^{\#}$, we get that $\mu^{\#}\trianglelefteq
\HH(\mu)$ and $h(\mu)\trianglelefteq\mu^{\#}$, hence
$\HH(\mu)=\mu^{\#}$, as required. This completes the proof of
Theorem \ref{main1}.

\begin{cor} For any $\lam\in\mathcal{P}_n$ and
$\mu\in\mathcal{K}_n$,
${[}\ts_{\bbQ}^{\lam}:{\td}_{\bbQ}^{\mu}{]}={[}\ts_{\bbQ}^{\widehat{\lam}}:{\td}_{\bbQ}^{\HH(\mu)}{]}$.
In particular,
${[}\ts_{\bbQ}^{\widehat{\mu}}:{\td}_{\bbQ}^{\HH(\mu)}{]}=1$ for
any $\mu\in\mathcal{K}_n$.
\end{cor}
\medskip\bigskip

\section{Proof of Theorem \ref{main2}}

In this section, we shall give the proof of Theorem \ref{main2}. Our
main tools are Dipper--Mathas's Morita equivalence results
(\cite{DM}) for Ariki--Koike algebras and their connections with
type $A$ affine Hecke algebras. Throughout this section, we keep the
same assumptions and notations as in Theorem \ref{main2}. That is,
$K=\mathbb{C}$,\, $q,\eps\in \mathbb{C}$ be such that $\eps$ is a
primitive $p$-th root of unity, $\eps^k=q^{\ell}$ is a primitive
$d$-th root of unity and $q$ is a primitive $(d\ell)$-th root of
unity, and $1<k<p$ is the smallest positive integer such that
$\eps^{k}\in\langle q\rangle$.\medskip

Let $\H_q(\BS_n)$ be the Iwahori--Hecke algebra associated to the
symmetric group $\BS_n$. Let $K[X_1^{\pm 1},\cdots,X_n^{\pm 1}]$ be
the ring of Laurent polynomials on $n$ indeterminates
$X_1,\cdots,X_n$.

\begin{dfn} The type $A$ affine Hecke algebra $\H_n^{\aff}$ is the $K$-algebra,
which as a $K$-linear space is isomorphic to $$
\H_q(\BS_n)\otimes_{K}K[X_1^{\pm 1},\cdots,X_n^{\pm 1}].
$$
The algebra structure is given by requiring that $\H_q(\BS_n)$ and
$K[X_1^{\pm 1},\cdots,X_n^{\pm 1}]$ are subalgebras and that
\addtocounter{equation}{1}\begin{equation}\label{equ51}
T_if-{\null}^{s_i}\negthickspace f T_i=(q-1)\frac{f-
{\null}^{s_i}\negthickspace f}{1-X_iX_{i+1}^{-1}}, \quad
\forall\,\,f\in K[X_1^{\pm 1},\cdots,X_n^{\pm 1}],
\end{equation}
Here $s_i\in\BS_n$ act on $K[X_1^{\pm 1},\cdots,X_n^{\pm 1}]$ by
permuting $X_i$ and $X_{i+1}$.
\end{dfn}

Note that the relation (\ref{equ51}) is equivalent to
$$\begin{aligned} & T_iX_{i}T_i=qX_{i+1},\quad \text{$\forall\,\, i$
with $1\leq i<n$,}\\ & T_iX_j=X_j T_i, \quad
\forall\,\,j\not\in\{i, i+1\},\end{aligned}
$$

Let $Q_1,\cdots,Q_p$ be elements of $K$. Let
$\H(p,n):=\H_{q,Q_1,\cdots,Q_p}(p,n)$ be the Ariki--Koike algebras
with parameters $\{q,Q_1,\cdots,Q_p\}$. It is well-known that there
is a surjective $K$-algebra homomorphism
$\varphi:\,\H_n^{\aff}\twoheadrightarrow\H(p,n)$ which is defined on
generators by
$$ T_i\mapsto T_i,\quad X_j\mapsto L_j,\,\,\forall\,1\leq
i<n,\,\,\forall\,1\leq j\leq n,
$$
where $L_j:=q^{1-j}T_{j-1}\cdots T_1T_0T_1\cdots T_{j-1}$ (the
$j$-th Murphy operator). As a consequence, every simple
$\H(p,n)$-module is a simple $\H_n^{\aff}$-module.\smallskip

We shall use Dipper--Mathas's explicit construction of Morita
equivalence for Ariki--Koike algebras. To this end, we need some
notations and definitions. Let $\{s_1,s_2,\cdots,s_{n-1}\}$ be the
set of basic transpositions in $\BS_n$. A word $w=s_{i_1}\cdots
s_{i_k}$ for $w\in\BS_{n}$ is a reduced expression if $k$ is
minimal; in this case we say that $w$ has length $k$ and we write
$\ell(w)=k$. Given a reduced expression $s_{i_1}\cdots s_{i_k}$ for
$w\in\BS_n$, we write $T_w=T_{i_1}\cdots T_{i_k}$, then $T_w$
depends only on $w$ and not on the choice of reduced expression. It
is well-known that $\H_{{q}}(\BS_n)$ is a free module with basis
$\{T_w|w\in\BS_n\}$. For each integer $0\leq a\leq n$, we define
$$
w_{n-a,a}=\underbrace{(s_{n-a}\cdots s_{n-1})}_{\text{$a$ times}}
\underbrace{(s_{n-a-1}\cdots s_{n-2})}_{\text{$a$ times}}\cdots
\underbrace{(s_{1}\cdots s_{a})}_{\text{$a$ times}}
$$ if
$a\not\in\{0,n\}$; or $w_{n-a,a}=1$ if $a\in\{0,n\}$.

Let $s$ be an integer with $1\leq s\leq p$ and such that $$
\prod_{1\leq i\leq s<j\leq p}\prod_{-n<a<n}(q^aQ_i-Q_j)\neq 0,
$$
in $K$.

\addtocounter{dfn}{1}
\begin{dfn} {\rm(\cite{DM})} For each integer $0\leq b\leq n$, let
$$\begin{aligned}
u_{n-b}^{-}:&=\prod_{t=1}^{s}(L_1-Q_t)(L_2-Q_t)\cdots
(L_{n-b}-Q_t),\\
u_{b}^{+}:&=\prod_{t=s+1}^{p}(L_1-Q_t)(L_2-Q_t)\cdots
(L_{b}-Q_t),\\
v_{b}:&=u_{n-b}^{-}T_{w_{n-b,b}}u_{b}^{+},\,\,V^b:=v_b\H(p,n).
\end{aligned} $$
\end{dfn}

Let $\H(s,b)=\H_{q,Q_1,\cdots,Q_s}(s,b)$ (resp.,
$\H(p-s,n-b)=\H_{q,Q_{s+1},\cdots,Q_p}(p-s,n-b)$) be the
Ariki--Koike algebra with parameters $\{q,Q_1,\cdots,Q_s\}$ and of
size $b$ (resp., with parameters $\{q,Q_{s+1},\cdots,Q_p\}$ and of
size $n-b$). Following \cite{DM}, let $\hat{\Pi}_{s,b}$ be the
functor from $\mmod_{\H(s,b)\otimes\H(p-s,n-b)}$ to
$\mmod_{\H(p,n)}$ given by $$
\hat{\Pi}_{s,b}(X):=X\otimes_{\H(s,b)\otimes\H(p-s,n-b)}V^{b},\,\,\,
\forall\,X\in \mmod_{\H(s,b)\otimes\H(p-s,n-b)}.$$ Note that here
the left $\bigl(\H(s,b)\otimes\H(p-s,n-b)\bigr)$-action on $V^{b}$
is well-defined because of the following very useful result from
\cite[(3.4)]{DM}.

\addtocounter{lem}{3}
\begin{lem}\label{lm54} {\rm(\cite[(3.4)]{DM})} Suppose that $0\leq b\leq n$.

(i) $$ T_i v_b=\begin{cases} v_b T_{i+b}, &\text{if $1\leq
i<n-b$,}\\
v_bT_{i-n+b}, &\text{if $n-b<i\leq n$.}\end{cases} $$

(ii) $$ L_k v_b=\begin{cases} v_b L_{k+b}, &\text{if $1\leq
k\leq n-b$,}\\
v_bL_{k-n+b}, &\text{if $n-b+1\leq k\leq n$.}\end{cases} $$
\end{lem}

Let $\bbQ_1:=(Q_1,\cdots,Q_s)$, $\bbQ_2:=(Q_{s+1},\cdots,Q_p)$. Let
$\bbQ:=(\bbQ_1,\bbQ_2)=(Q_1,\cdots,Q_p)$. Suppose that
$D^{\lam^{[1]}}_{\bbQ_1}\neq 0$ (resp., $D^{\lam^{[2]}}_{\bbQ_2}\neq
0$) is an irreducible $\H(s,b)$-module (resp.,
$\H(p-s,n-b)$-module), where $\lam^{[1]}$ (resp., $\lam^{[2]}$) is
an $s$-multipartition of $b$ (resp., $(p-s)$-multipartition of
$n-b$). Let $\lam:=(\lam^{[1]},\lam^{[2]})$ (concatenation of
ordered tuples), which is a $p$-multipartition of $n$. By \cite{DM},
$D^{\lam}_{\bbQ}\neq 0$ is an irreducible $\H(p,n)$-module, and
$$
\hat{\Pi}_{s,b}\Bigl(D^{\lam^{[1]}}_{\bbQ_1}\otimes
D^{\lam^{[2]}}_{\bbQ_2}\Bigr)=\Bigl(D^{\lam^{[1]}}_{\bbQ_1}\otimes
D^{\lam^{[2]}}_{\bbQ_2}\Bigr)\otimes_{\H(s,b)\otimes\H(p-s,n-b)}V^b\cong
D^{\lam}_{\bbQ}.
$$

Let $\H_b^{\aff}$ (resp., $\H_{n-b}^{\aff}$) be the standard
parabolic subalgebra of $\H_n^{\aff}$ generated by
$T_1,\cdots,T_{b-1},X_1,\cdots,X_b$ (resp., by
$T_{b+1},\cdots,T_{n-1},X_{b+1},\cdots,X_n$). Then
$D^{\lam^{[1]}}_{\bbQ_1}$ (resp., $D^{\lam^{[2]}}_{\bbQ_2}$)
naturally becomes an irreducible $\H_{b}^{\aff}$-module (resp.,
$\H_{n-b}^{\aff}$-module). We have that

\addtocounter{prop}{4}
\begin{prop} \label{thm55} There is an $\H_n^{\aff}$-module isomorphism $$
D^{\lam}_{\bbQ}\cong\Ind_{\H_{b}^{\aff}\otimes\H_{n-b}^{\aff}}^{\H_n^{\aff}}\Bigl(
D^{\lam^{[1]}}_{\bbQ_1}\otimes D^{\lam^{[2]}}_{\bbQ_2}\Bigr).
$$
\end{prop}

\noindent {Proof:} \,By our previous discussion, it suffices to
show that $$ \begin{aligned}&\Bigl(D^{\lam^{[1]}}_{\bbQ_1}\otimes
D^{\lam^{[2]}}_{\bbQ_2}\Bigr)\otimes_{\H_{b}^{\aff}\otimes\H_{n-b}^{\aff}}\H_{n}^{\aff}
\cong\\
&\qquad\qquad\Bigl(D^{\lam^{[1]}}_{\bbQ_1}\otimes
D^{\lam^{[2]}}_{\bbQ_2}\Bigr)\otimes_{\H(s,b)\otimes\H(p-s,n-b)}v_b\H(p,n),
\end{aligned}$$ where $v_b\H(p,n)$ is regarded as right
$\H_{n}^{\aff}$-module via the natural surjective homomorphism
$\varphi:\,\H_{n}^{\aff}\twoheadrightarrow\H(p,n)$.

In fact, by Lemma \ref{lm54}, it is easy to see that the following
map $$ \bigl(x\otimes
y\bigr)\otimes_{\H_{b}^{\aff}\otimes\H_{n-b}^{\aff}} h\mapsto
\bigl(x\otimes y\bigr)\otimes_{\H(s,b)\otimes\H(p-s,n-b)} v_b h
$$
extends naturally to a well-defined surjective right
$\H_{n}^{\aff}$-module homomorphism. Now comparing their
dimensions (see \cite[(4.8)]{DM}), we proves the theorem.\hfill\qed\medskip

Now we suppose that
$\bQ=\bQ_1\sqcup\bQ_2\sqcup\cdots\sqcup\bQ_{\kappa}$ (disjoint
union) such that $Q_i, Q_j$ are in the same $q$-orbit only if
$Q_i, Q_j\in\bQ_c$ for some $1\leq c\leq\kappa$. For each integer
$i$ with $1\leq i\leq\kappa$, let $D^{\lam^{[i]}}_{\bbQ_i}\neq 0$
be an irreducible $\H(p_i,b_i)$-module, where $p_i=|\bbQ_i|$,
$\lam^{[i]}$ is a $p_i$-multipartition of $b_i$,
$\sum_{i=1}^{\kappa}b_i=n$. Let
$\lam:=(\lam^{[1]},\cdots,\lam^{[\kappa]})$ (concatenation of
ordered tuples), which is a $p$-multipartition of $n$.

\addtocounter{cor}{5}
\begin{cor}\label{cor56} With the above assumptions and notations, we have an $\H_n^{\aff}$-module isomorphism $$
D^{\lam}_{\bbQ}\cong\Ind_{\H_{b_1}^{\aff}\otimes\cdots\otimes\H_{b_{\kappa}}^{\aff}}^{\H_n^{\aff}}\Bigl(
D^{\lam^{[1]}}_{\bbQ_1}\otimes\cdots\otimes
D^{\lam^{[\kappa]}}_{\bbQ_{\kappa}}\Bigr).
$$
\end{cor}

\noindent {Proof:} \,This follows from Proposition \ref{thm55} and the
associativity of tensor product induction functor.\hfill\qed\smallskip

Now we return to our setup in Theorem \ref{main2}. For each $1\leq
i\leq k$, we set $\bbQ_i=(\veps^{i-1},\veps^{k+i-1},\cdots,
\veps^{(d-1)k+i-1})$. $\bQ=\bQ_1\sqcup\cdots\sqcup\bQ_k$ is a
partition of the parameters set $\bQ$ into different $q$-orbits. Let
$\bbQ=\bigl(\bbQ_1,\bbQ_2,\cdots,\bbQ_{k}\bigr)$ (concatenation of
ordered tuples). For each $p$-multipartition
$\lam=(\lam^{(1)},\cdots,\lam^{(p)})$ of $n$, we write
$$
\lam^{[i]}=(\lam^{((i-1)d+1)},\lam^{((i-1)d+2)},\cdots,\lam^{(id)}),\,\,\,\text{for}\,\,i=1,2,\cdots,k.
$$ and recall $\theta$ is the map $\lam\mapsto
(\lam^{[1]},\cdots,\lam^{[k]})$. Let $n_i=|\lam^{[i]}|$ for each
$1\leq i\leq k$.

\bigskip
\noindent {\bf Proof of Theorem \ref{main2}:}\quad We define $$
\bbQ_{1}^{\ast}:=(\eps^k,\eps^{2k},\cdots,\eps^{(d-1)k},1).
$$
By Lemma \ref{lm24}, we have that $$
\Bigl(\td_{\bbQ}^{\lam}\Bigr)^{\sigma}\cong\td_{(\bbQ_{2},\bbQ_{3},\cdots,\bbQ_{k},\bbQ_{1}^{\ast})}^{\lam}.
$$
Applying Corollary \ref{cor56}, we get that $$\begin{aligned}
\Bigl(\td_{\bbQ}^{\lam}\Bigr)^{\sigma}&\cong\td_{(\bbQ_{2},\bbQ_{3},\cdots,\bbQ_{k},\bbQ_{1}^{\ast})}^{\lam}\\
&\cong
\Ind_{\H_{n_1}^{\aff}\otimes\cdots\otimes\H_{n_{k-1}}^{\aff}\otimes\H_{n_{k}}^{\aff}}^{\H_n^{\aff}}\Bigl(
D^{\lam^{[1]}}_{\bbQ_2}\otimes
D^{\lam^{[2]}}_{\bbQ_3}\otimes\cdots\otimes
D^{\lam^{[k-1]}}_{\bbQ_{k}}\otimes
D^{\lam^{[k]}}_{\bbQ_1^{\ast}}\Bigr).\end{aligned}
$$
By \cite[(5.12)]{V}, the righthand side module has the same
composition factors as $$
\Ind_{\H_{n_k}^{\aff}\otimes\H_{n_{1}}^{\aff}\otimes\cdots\otimes\H_{n_{k-1}}^{\aff}}^{\H_n^{\aff}}\Bigl(
D^{\lam^{[k]}}_{\bbQ_1^{\ast}}\otimes
D^{\lam^{[1]}}_{\bbQ_2}\otimes
D^{\lam^{[2]}}_{\bbQ_3}\otimes\cdots\otimes
D^{\lam^{[k-1]}}_{\bbQ_{k}}\Bigr).
$$
In particular, as both modules are irreducible, these two modules
are in fact isomorphic to each other. Therefore, $$
\Bigl(\td_{\bbQ}^{\lam}\Bigr)^{\sigma}\cong\Ind_{\H_{n_k}^{\aff}\otimes\H_{n_{1}}^{\aff}\otimes\cdots\otimes\H_{n_{k-1}}^{\aff}}^{\H_n^{\aff}}\Bigl(
D^{\lam^{[k]}}_{\bbQ_1^{\ast}}\otimes
D^{\lam^{[1]}}_{\bbQ_2}\otimes
D^{\lam^{[2]}}_{\bbQ_3}\otimes\cdots\otimes
D^{\lam^{[k-1]}}_{\bbQ_{k}}\Bigr).
$$

Now again by Lemma \ref{lm24}, $$
D^{\lam^{[k]}}_{\bbQ_1^{\ast}}\cong
\Bigl(\td_{\bbQ_1}^{\lam^{[k]}}\Bigr)^{\sigma'}, $$ where $\sigma'$
denotes the $K$-algebra automorphism of the Ariki--Koike algebra
$$\H_{q}(d,n_k):=\H_{q,1,\eps',\cdots,(\eps')^{d-1}}(d,n_k)$$
(where $\eps'=\eps^k$) which is defined on generators by
$\sig(T_0)=\eps' T_0, \sig(T_i)=T_i$, for $i=1,2,\cdots,n_k-1$. By
Theorem \ref{main1} and Corollary \ref{maincor}, we deduce that
$$
D^{\lam^{[k]}}_{\bbQ_1^{\ast}}\cong
\Bigl(\td_{\bbQ_1}^{\lam^{[k]}}\Bigr)^{\sigma'}\cong
D^{\HH'(\lam^{[k]})}_{\bbQ_1},
$$
where $\HH'$ is as defined in Corollary \ref{maincor}. Therefore
$$\begin{aligned}
\Bigl(\td_{\bbQ}^{\lam}\Bigr)^{\sigma}&\cong\Ind_{\H_{n_k}^{\aff}\otimes\H_{n_{1}}^{\aff}\otimes\cdots\otimes\H_{n_{k-1}}^{\aff}}^{\H_n^{\aff}}\Bigl(
D^{\lam^{[k]}}_{\bbQ_1^{\ast}}\otimes
D^{\lam^{[1]}}_{\bbQ_2}\otimes
D^{\lam^{[2]}}_{\bbQ_3}\otimes\cdots\otimes
D^{\lam^{[k-1]}}_{\bbQ_{k}}\Bigr)\\
&\cong\Ind_{\H_{n_k}^{\aff}\otimes\H_{n_{1}}^{\aff}\otimes\cdots\otimes\H_{n_{k-1}}^{\aff}}^{\H_n^{\aff}}\Bigl(
D^{\HH'(\lam^{[k]})}_{\bbQ_1}\otimes
D^{\lam^{[1]}}_{\bbQ_2}\otimes
D^{\lam^{[2]}}_{\bbQ_3}\otimes\cdots\otimes
D^{\lam^{[k-1]}}_{\bbQ_{k}}\Bigr)\\
&\cong
\td_{\bbQ}^{(\HH'(\lam^{[k]}),\lam^{[1]},\cdots,\lam^{[k-1]})}.
\end{aligned}
$$
It follows that $$
\theta\Bigl(\HH(\lam)\Bigr)=\Bigl(\HH'(\lam^{[k]}),\lam^{[1]},\cdots,\lam^{[k-1]}\Bigr),
$$
as required. This completes the proof of Theorem \ref{main2}.\hfill\qed
\medskip\bigskip

\section{Closed formula for the number of simple $\H_q(p,p,n)$-modules}

The purpose of this section is to give the second two main results
(Theorem \ref{mainthm3} and Theorem \ref{mainthm4}) of this paper,
which yield explicit formula for the number of simple modules over
the cyclotomic Hecke algebra of type $G(p,p,n)$ in the non-separated
case. Note that in the separated case one can easily write down an
explicit formula by using the result \cite[(5.7)]{Hu3}. In the
non-separated case we shall apply Theorem \ref{main1} and Theorem
\ref{main2} as well as Naito--Sagaki's work (\cite{NS1},\cite{NS2})
on Lakshmibai--Seshadri paths fixed by diagram
automorphisms.\medskip

Recall our definitions of ${K}_n$ and $\HH$ in the second paragraph
of Section 3.

\begin{dfn} For each $\lam\in{K}_n$, let $$
o_{\HH}(\lam):=\min\bigl\{1\leq m\leq
p\bigm|\HH^m(\lam)=\lam\bigr\}.
$$ For each integer $m$ with $1\leq m\leq p$, we define
$$\begin{aligned}
&\widetilde{\Sigma}(m):=\bigl\{\lam\in{K}_n\bigm|\HH^{m}(\lam)=\lam\bigr\},\,\,\,
\widetilde{N}(m):=\#\widetilde{\Sigma}(m),\\
& N(m):=\#\bigl\{\lam\in{K}_n\bigm|o_{\HH}(\lam)=m\bigr\}.
\end{aligned}
$$
\end{dfn}

We use the notation $\#\Irr\bigl(\H_q(p,p,n)\bigr)$ (resp.,
$\#\Irr\bigl(\H_{q}(p,n)\bigr)$) to denote the number of simple
$\H_q(p,p,n)$-modules (resp., simple $\H_{q}(p,n)$-modules). By
Lemma \ref{thm31}, we know that \addtocounter{equation}{1}
\begin{equation}\label{equa62}\begin{split}
\#\Irr\bigl(\H_q(p,p,n)\bigr)&=\frac{1}{p}\Bigl\{\#\Irr\bigl(\H_{q}(p,n)\bigr)
-\sum_{1\leq m<p,
m|p}N(m)\Bigr\}\\
&\qquad+\sum_{1\leq
m<p,m|p}\frac{N(m)}{m}\frac{p}{m}.\end{split}\end{equation} Note
that by \cite{AM}, the number $\#\Irr\bigl(\H_{q}(p,n)\bigr)$ is
explicitly known. Therefore, to get a formula for
$\#\Irr\bigl(\H_q(p,p,n)\bigr)$, it suffices to derive a formula
for ${N}(\widetilde{m})$ for each integer $1\leq \widetilde{m}<p$
satisfying $\widetilde{m}|p$.

Let $\mu$ be the M\"{o}bius function
$\mu:\,\mathbb{N}\rightarrow\{0,1,-1\}$ which is given by $$
\mu(a)=\begin{cases} 1, &\text{if $a=1$,}\\
(-1)^{s}, &\begin{matrix}\text{if $a=p_1\cdots p_s$, where
$\{p_i\}_{1\leq i\leq s}$ are}\\[-4pt] \text{pairwise
different prime numbers,}\end{matrix}\\
0 &\text{otherwise}
\end{cases}
$$
Since $\widetilde{N}(m)=\sum_{1\leq a\leq m, a|m}N(a)$, it follows
from M\"{o}bius inversion formula that
\begin{equation}\label{equa63}
N(\widetilde{m})=\sum_{1\leq m\leq \widetilde{m},
m|\widetilde{m}}\mu(\widetilde{m}/m)\widetilde{N}(m).\end{equation}
Therefore, it suffices to derive a formula for $\widetilde{N}(m)$
for each integer $1\leq m<p$ satisfying $m|p$. To this end, we have
to use Naito--Sagaki's work (\cite{NS1},\cite{NS2}) on
Lakshmibai--Seshadri paths fixed by diagram automorphisms.

For the moment, we assume the following setup. That is,
$K=\mathbb{C}$,\, $q,\eps\in \mathbb{C}$ be such that
$\eps=q^{\ell}$ is a primitive $p$-th root of unity, and $q$ is a
primitive $p\ell$-th root of unity. We identify $K_n$ with
$\mathcal{K}_n$, the set of Kleshchev $p$-multipartitions with
respect to $\{q,1,\eps,\cdots,\eps^{p-1}\}$. Let $\mathfrak{g}$ be
the Kac--Moody algebra over $\mathbb{C}$ associated to a
symmetrizable generalized Cartan matrix $(a_{i,j})_{i,j\in I}$ of
finite size. Let $\mathfrak{h}$ be its Cartan subalgebra, and
$W\subset\GL(\mathfrak{h}^{\ast})$ be its Weyl group. Let
$\{\alpha_i^{\vee}\}_{0\leq i\leq n-1}$ be the set of simple coroots
in $\mathfrak{h}$. Let
$\mathcal{X}:=\bigl\{\Lambda\in\mathfrak{h}^{\ast}\bigm|
\Lambda(\alpha_i^{\vee})\in\Z,\,\forall\,0\leq i<n\bigr\}$ be the
weight lattice. Let
$\mathcal{X}^{+}:=\bigl\{\Lambda\in\mathcal{X}\bigm|\Lambda(\alpha_i^{\vee})\geq
0,\,\forall\,0\leq i<n\bigr\}$ be the set of integral dominant
weights. Let
$\mathcal{X}_{\mathbb{R}}:=\mathcal{X}\otimes_{\Z}\mathbb{R}$, where
$\mathbb{R}$ is the real number field. Assume that
$\Lambda\in\mathcal{X}^{+}$. P. Littelmann introduced (\cite{Li1},
\cite{Li2}) the notion of Lakshmibai--Seshadri paths (L-S paths for
short) of class $\Lambda$, which are piecewise linear, continuous
maps $\pi:[0,1]\rightarrow\mathcal{X}_{\mathbb{R}}$ parameterized by
pairs $(\underline{\nu},\underline{a})$ of a sequence
$\underline{\nu}: \nu_1>\nu_2>\cdots>\nu_s$ of elements of
$W\Lambda$, where $>$ is the ``relative Bruhat order" (see
\cite[Section 4]{Li2}) on $W\Lambda$, and a sequence $\underline{a}:
0=a_0<a_1<\cdots<a_s=1$ of rational numbers with a certain
condition, called the chain condition. The set $\mathbb{B}(\Lambda)$
of all L-S paths of class $\Lambda$ is called the path model for the
irreducible integrable highest weight module $L(\Lambda)$ of highest
weight $\Lambda$ over $\mathfrak{g}$. It is a remarkable fact that
$\mathbb{B}(\Lambda)$ has a canonical crystal structure which is
isomorphic to the crystal associated to the irreducible integrable
highest weight module of highest weight $\Lambda$ over the quantum
affine algebra $U'_v(\mathfrak{g})$.

Now let $\mathfrak{g}=\ksl_{p\l}$, the affine Kac--Moody algebra of
type $A_{p\l-1}^{(1)}$. The generalized Cartan matrix
$(a_{i,j})_{i,j\in I}$ of $\mathfrak{g}$ was indexed by the finite
set $I:=\Z/p\l\Z$. Let $1\leq m<p$ be an integer satisfying $m|p$.
Let $\omega:\,I\rightarrow I$ be an automorphism of order $p/m$
defined by $\bar{i}=i+p\l\Z\mapsto
\bar{i}-\overline{m\l}=i-m\l+p\l\Z$ for any $\bar{i}\in I$. Clearly
$\omega$ is a Dynkin diagram automorphism in the sense of
\cite[\S1.2]{NS1} (i.e., satisfying
$a_{\omega(i),\omega(j)}=a_{i,j}$, $\forall\,i,j\in I$). By
\cite{FSS}, $\omega$ induces a Lie algebra automorphism (which is
called diagram outer automorphism) $\omega\in\Aut(\mathfrak{g})$ of
order $p/m$ and a linear automorphism
$\omega^{\ast}\in\GL(\mathfrak{h}^{\ast})$ of order $p/m$.

Following \cite{FRS} and \cite[\S1.3]{NS1}, we set
$c_{i,j}:=\sum\limits_{t=0}^{N_j-1}a_{i,\omega^t(j)}$, where
$N_j:=\#\bigl\{\omega^t(i)\bigm|t\geq 0\bigr\}$, $i,j\in I$. We
choose a complete set $\widehat{I}$ of representatives of the
$\omega$-orbits in $I$, and set
$\check{I}:=\bigl\{i\in\widehat{I}\bigm|c_{i,i}>0\bigr\}$. We put
$\hat{a}_{i,j}:=2c_{i,j}/c_{j}$ for $i,j\in\widehat{I}$, where
$c_i:=c_{ii}$ if $i\in\check{I}$, and $c_i:=2$ otherwise. Then
$(\hat{a}_{i,j})_{i,j\in\widehat{I}}$ is a symmetrizable
Borcherds--Cartan matrix in the sense of \cite{Bo}, and (if
$\check{I}\neq\emptyset$) its submatrix
$(\hat{a}_{i,j})_{i,j\in\check{I}}$ is a generalized Cartan matrix
of affine type. Let $\widehat{\mathfrak{g}}$ be the generalized
Kac--Moody  algebra over $\mathbb C$ associated to
$(\hat{a}_{i,j})_{i,j\in\widehat{I}}$, with Cartan subalgebra
$\widehat{\mathfrak{h}}$, Chevalley generators
$\{\hat{x}_i,\hat{y}_i\}_{i\in\widehat{I}}$. The orbit Lie algebra
$\check{\mathfrak{g}}$ is defined to be the subalgebra of
$\widehat{\mathfrak{g}}$ generated by $\widehat{\mathfrak{h}}$ and
$\hat{x}_i,\hat{y}_i$ for $i\in\check{I}$, which is a usual
Kac--Moody algebra.

\addtocounter{lem}{3}
\begin{lem} With the above assumptions and notations, we have that
$$
\check{\mathfrak{g}}=\begin{cases} \ksl_{m\l}, &\text{if $m\l>1$,}\\
\mathbb{C}, &\text{if $m=\l=1$.}
\end{cases}
$$
\end{lem}

\noindent {Proof:} \,This follows from direct
verification.\hfill\qed\medskip

We define
$\bigl(\mathfrak{h}^{\ast}\bigr)^{\circ}:=\bigl\{\Lambda\in\mathfrak{h}^{\ast}
\bigm|\omega^{\ast}(\Lambda)=\Lambda\bigr\}$.
$\widetilde{W}:=\bigl\{w\in
W\bigm|\omega^{\ast}w=w\omega^{\ast}\bigr\}$. To distinguish with
the objects for $\mathfrak{g}$, the objects for the obit Lie algebra
$\check{\mathfrak{g}}$ will always have the symbol ``$\vee{}$'' on
the head. For example, $\check{\mathfrak{h}}$ denotes the Cartan
subalgebra of $\check{\mathfrak{g}}$, $\check{W}$ the Weyl group of
$\check{\mathfrak{g}}$, $\{\check{\Lambda}_i\}_{0\leq i\leq m\l-1}$
the set of fundamental dominant weights in
$\check{\mathfrak{h}}^{\ast}$. There exist a linear automorphism
$P_{\omega}^{\ast}:\,\check{\mathfrak{h}}^{\ast}\rightarrow
\bigl(\mathfrak{h}^{\ast}\bigr)^{\circ}$ and a group isomorphism
$\Theta:\,\check{W}\rightarrow\widetilde{W}$ such that
$\Theta(\check{w})=P_{\omega}^{\ast}\check{w}\bigl(P_{\omega}^{\ast}\bigr)^{-1}$
for each $w\in\check{W}$. By \cite[\S6.5]{FSS}, for each $0\leq
i<m\l$, $$
P_{\omega}^{\ast}(\check{\Lambda}_i)=\Lambda_{i}+\Lambda_{i+m\l}+\Lambda_{i+2m\l}+\cdots+\Lambda_{i+(p-m)\l}+C\delta,
$$
where $C\in\mathbb{Q}$ is some constant depending on $\omega$,
$\delta$ denotes the null root of $\mathfrak{g}$. Let
$\check{\Lambda}=\check{\Lambda}_0+\check{\Lambda}_{\l}+\check{\Lambda}_{2\l}+\cdots+\check{\Lambda}_{(m-1)\l}$.
Let
$\Lambda:=\Lambda_0+\Lambda_{\l}+\Lambda_{2\l}+\cdots+\Lambda_{(p-1)\l}$.
Then it follows that
$P_{\omega}^{\ast}(\check{\Lambda})=\Lambda+C'\delta$, for some
$C'\in\mathbb{Q}$.  \smallskip

Let $\mathbb{B}(\Lambda)$ (resp.,
$\mathbb{B}\bigl(P_{\omega}^{\ast}(\check{\Lambda})\bigr)$) be the
set of all L-S paths of class $\Lambda$ (resp., of class
$P_{\omega}^{\ast}(\check{\Lambda})$). Let $\pi_{\Lambda}$ (resp.,
$\pi_{P_{\omega}^{\ast}(\check{\Lambda})}$) be the straight path
joining $0$ and $\Lambda$ (resp., $0$ and
$P_{\omega}^{\ast}(\check{\Lambda})$). Let $\widetilde{E}_i,
\widetilde{F}_i$ denote the raising root operator and the lowering
root operator (see \cite{Li1} and \cite{Li2}) with respect to the simple root $\alpha_i$.

\begin{lem} The map which sends $\pi_{P_{\omega}^{\ast}(\check{\Lambda})}$
to $\pi_{\Lambda}$ extends to a bijection $\beta$ from
$\mathbb{B}\bigl(P_{\omega}^{\ast}(\check{\Lambda})\bigr)$ onto
$\mathbb{B}(\Lambda)$ such that $$
\beta\bigl(\widetilde{F}_{i_1}\cdots
\widetilde{F}_{i_s}\pi_{P_{\omega}^{\ast}(\check{\Lambda})}\bigr)=\widetilde{F}_{i_1}\cdots
\widetilde{F}_{i_s}\pi_{\Lambda},
$$
for any $i_1,\cdots,i_s\in\Z/{p\l\Z}$.
\end{lem}

\noindent {Proof:} \,This follows from the fact that
$P_{\omega}^{\ast}(\check{\Lambda})-\Lambda\in\mathbb{Q}\delta$
and the definitions of
$\mathbb{B}\bigl(P_{\omega}^{\ast}(\check{\Lambda})\bigr)$ and
$\mathbb{B}(\Lambda)$.\hfill\qed\medskip

Henceforth we shall identify
$\mathbb{B}\bigl(P_{\omega}^{\ast}(\check{\Lambda})\bigr)$ with
$\mathbb{B}(\Lambda)$. The action of $\omega^{\ast}$ on
$\mathfrak{h}^{\ast}$ naturally extends to the set
$\mathbb{B}\bigl(P_{\omega}^{\ast}(\check{\Lambda})\bigr)$ (and
hence to the set $\mathbb{B}\bigl(\Lambda\bigr)$). By
\cite[(3.1.1)]{NS2}, if
$\widetilde{F}_{i_1}\widetilde{F}_{i_2}\cdots \widetilde{F}_{i_s}
\pi_{{\Lambda}}\in {\mathbb{B}}({\Lambda})$, then
\addtocounter{equation}{2}
\begin{equation}\label{equa66}
\omega^{\ast}\bigl(\widetilde{F}_{i_1}\widetilde{F}_{i_2}\cdots
\widetilde{F}_{i_s}
\pi_{{\Lambda}}\bigr)=\widetilde{F}_{i_1+m\l}\widetilde{F}_{i_2+m\l}\cdots
\widetilde{F}_{i_s+m\l} \pi_{{\Lambda}}.
\end{equation}

We denote by $\mathbb{B}^{\circ}\bigl(\Lambda\bigr)$ the set of
all L-S paths of class $\Lambda$ that are fixed by
$\omega^{\ast}$. For $\check{\mathfrak{g}}$, we denote by
$\widetilde{e}_i, \widetilde{f}_i$ the raising root operator and
the lowering root operator with respect to the simple root
$\alpha_i$. Let $\pi_{\check{\Lambda}}$ be the straight path
joining $0$ and $\check{\Lambda}$. By \cite[(4.2)]{NS1}, the
linear map $P_{\omega}^{\ast}$ naturally extends to a map from
$\check{\mathbb{B}}(\check{\Lambda})$ to
$\mathbb{B}^{\circ}\bigl(\Lambda\bigr)$ such that if
$\widetilde{f}_{i_1}\widetilde{f}_{i_2}\cdots \widetilde{f}_{i_s}
\pi_{\check{\Lambda}}\in \check{\mathbb{B}}(\check{\Lambda})$,
then
$$\begin{aligned}
&P_{\omega}^{\ast}\bigl(\widetilde{f}_{i_1}\widetilde{f}_{i_2}\cdots
\widetilde{f}_{i_s}
\pi_{\check{\Lambda}}\bigr)=\widetilde{F}_{i_1}\widetilde{F}_{i_1+m\l}\cdots\widetilde{F}_{i_1+(p-m)\l}
\widetilde{F}_{i_2}\widetilde{F}_{i_2+m\l}\cdots\widetilde{F}_{i_2+(p-m)\l}\cdots\\
&\qquad\qquad\qquad\qquad\qquad\widetilde{F}_{i_s}\widetilde{F}_{i_s+m\l}\cdots\widetilde{F}_{i_s+(p-m)\l}
\pi_{\Lambda}.\end{aligned}
$$

\addtocounter{lem}{1}
\begin{lem} \label{thm67} {\rm(\cite[(4.2),(4.3)]{NS1})}
$\mathbb{B}^{\circ}\bigl(\Lambda\bigr)=P_{\omega}^{\ast}\bigl(
\check{\mathbb{B}}(\check{\Lambda})\bigr)$.
\end{lem}

Note that both $\check{\mathbb{B}}(\check{\Lambda})$ and
$\mathbb{B}\bigl(\Lambda\bigr)$ have a canonical crystal structure
with the raising and lowering root operators playing the role of
Kashiwara operators. They are isomorphic to the crystals associated
to the irreducible integrable highest weight modules
$\check{L}(\check\Lambda)$ of highest weight $\check\Lambda$
 over $U'_v(\check{\mathfrak{g}})$ and $L\bigl(\Lambda\bigr)$
of highest weight $\Lambda$ over $U'_v(\mathfrak{g})$
 respectively. Henceforth, we identify
them without further comments. Let $v_{\check{\Lambda}}$ (resp.,
$v_{\Lambda}$) denotes the unique highest weight vector of highest
weight $\check{\Lambda}$ (resp., of highest weight $\Lambda$) in
$\check{\mathbb{B}}(\check{\Lambda})$ (resp., in
$\mathbb{B}(\Lambda)$). Therefore, by (\ref{equa66}) and Lemma
\ref{thm67}, we get that

\addtocounter{cor}{7}
\begin{cor} With the above assumptions and notations, there is an
injection $\eta$ from the set
$\check{\mathbb{B}}(\check{\Lambda})$ of crystal bases to the set
$\mathbb{B}(\Lambda)$ of crystal bases such that $$\begin{aligned}
&\eta\bigl(\widetilde{f}_{i_1}\widetilde{f}_{i_2}\cdots
\widetilde{f}_{i_s}
v_{\check{\Lambda}}\bigr)\equiv\widetilde{F}_{i_1}\widetilde{F}_{i_1+m\l}\cdots\widetilde{F}_{i_1+(p-m)\l}
\widetilde{F}_{i_2}\widetilde{F}_{i_2+m\l}\cdots\widetilde{F}_{i_2+(p-m)\l}\cdots\\
&\qquad\qquad\qquad\qquad\qquad\widetilde{F}_{i_s}\widetilde{F}_{i_s+m\l}\cdots\widetilde{F}_{i_s+(p-m)\l}
v_{\Lambda}\pmod{{v L(\Lambda)_{A}}},\end{aligned}
$$
and the image of $\eta$ consists of all crystal basis element
$\widetilde{F}_{i_1}\cdots \widetilde{F}_{i_t}v_{\Lambda}+v
L(\Lambda)_{A}$ satisfying $\widetilde{F}_{i_1}\cdots
\widetilde{F}_{i_t}v_{\Lambda}\equiv \widetilde{F}_{i_1+m\l}\cdots
\widetilde{F}_{i_t+m\l}v_{\Lambda}\pmod{{v L(\Lambda)_{A}}}. $
\end{cor}

We translate the language of crystal bases into the language of
Kleshchev multipartitions, we get the following combinatorial
result, which seems of independent interest.

\begin{cor} \label{cor609} Let $\HH$ be as in Theorem \ref{main1}. Let
$\eps':=\eps^{p/m}$. Let $q'=\sqrt[\l]{\eps'}$, which is a primitive $m\l$-th root of unity.
Then there exists a bijection
$\eta:\,\check{\lam}\mapsto\lam$ from the set of Kleshchev
$m$-multipartitions $\check{\lam}$ of $nm/p$ with respect to
$(q',1,\eps',\eps'^{2},\cdots,(\eps')^{m-1})$ onto the set of
Kleshchev $p$-multipartitions $\lam$ of $n$ with respect to
$(q,1,\eps,\eps^{2},\cdots,\eps^{p-1})$ satisfying
$\HH^m(\lam)=\lam$, such that if
$$
\underbrace{(\emptyset,\cdots,\emptyset)}_{m}\overset{r_1}{\twoheadrightarrow}\cdot
\overset{r_2}{\twoheadrightarrow}\cdot \cdots\cdots
\overset{r_s}{\twoheadrightarrow}\check\lam $$ is a path from
$\underbrace{(\emptyset,\cdots,\emptyset)}_{m}$ to $\check\lam$ in
Kleshchev's good lattice with respect to
$(q',1,\eps',(\eps')^{2},\cdots,(\eps')^{m-1})$, where $s:=nm/p$, then the sequence
$$\begin{aligned}&\underbrace{(\emptyset,\cdots,\emptyset)}_{p}
\overset{r_1}{\twoheadrightarrow}\cdot\overset{m\l+r_1}{\twoheadrightarrow}\cdot
\cdots\overset{(p-m)\l+r_1}{\twoheadrightarrow}\cdot
\overset{r_2}{\twoheadrightarrow}\cdot
\overset{m\l+r_2}{\twoheadrightarrow}\cdot
\cdots\overset{(p-m)\l+r_2}{\twoheadrightarrow}\cdot\\
&\qquad \qquad\cdots\,\cdot
\overset{r_s}{\twoheadrightarrow}\cdot\overset{m\l+r_s}{\twoheadrightarrow}\cdots
\overset{(p-m)\l+r_s}{\twoheadrightarrow}\lam\end{aligned}
$$ defines a path in Kleshchev's good lattice (w.r.t.,
$(q,1,\eps,\eps^2,\cdots,\eps^{p-1})$) satisfying
$\HH^m(\lam)=\lam$.
\end{cor}

\noindent {Proof:} \,This follows from (\ref{equa63}), Theorem
\ref{main1} and the realization of crystal graph as Kleshchev's good
lattice.\hfill\qed\medskip

We remark that one can derive a similar combinatorial result for FLOTW $p$-partitions
by using the same arguments.

\addtocounter{thm}{9}
\begin{thm} \label{mainthm3} Suppose that $K=\mathbb{C}$,\, $q,\eps,q',\eps'\in \mathbb{C}$ be such that
$\eps=q^{\ell}$ (resp., $\eps'=(q')^{\ell}$) is a primitive $p$-th
(resp., primitive $m$-th) root of unity, and $q$ (resp., $q'$) is a
primitive $p\ell$-th (resp., primitive $m\l$-th) root of unity. Let
$1\leq m\leq p$ be an integer such that $m|p$. Then
$\widetilde{N}(m)$ is equal to the number of simple
$\H_{q',1,\eps',\cdots,(\eps')^{m-1}}(m,mn/p)$-modules, where
$\H_{q',1,\eps',\cdots,(\eps')^{m-1}}(m,mn/p)$ is the Ariki--Koike
algebra with parameters $(q',1,\eps',\cdots,(\eps')^{m-1})$ and of
size $mn/p$. In particular, for each integer $1\leq \widetilde{m}<p$
such that $\widetilde{m}|p$, we get that \addtocounter{equation}{4}
\begin{equation}\label{equa611}\begin{split}
N(\widetilde{m})&=\sum_{1\leq
m\leq\widetilde{m}, m|\widetilde{m}}\mu(\widetilde{m}/m)\widetilde{N}(m)\\
&=\sum_{1\leq m\leq\widetilde{m},
m|\widetilde{m}}\mu(\widetilde{m}/m)\Bigl(\#\Irr\bigl(\H_{q',1,\eps',\cdots,(\eps')^{m-1}}(m,mn/p)\bigr)\Bigr).
\end{split}\end{equation}
\end{thm}

In the special case where $\l=m=1$, the set
$\check{\mathbb{B}}(\check{\Lambda})$ contains only one
element---a highest weight vector, which corresponds to the empty
partition $\emptyset$. Hence we have

\addtocounter{cor}{2}
\begin{cor} If $\l=1$, then for any integer $n\geq 1$, $$
\bigl\{\lam\in{K}_n\bigm|\HH(\lam)=\lam\bigr\}=\emptyset.
$$
Equivalently, $\widetilde{N}(1)=0=N(1)$.
\end{cor}

\begin{cor} With the same assumption as in Theorem \ref{mainthm3}, if
$(p,n)=1$, then $\widetilde{N}(m)=0$ for any integer $1\leq m<p$.
In this case, for any irreducible $\H_q(p,n)$-module $D$,
$D\downarrow_{\H_q(p,p,n)}$ remains irreducible.
\end{cor}

\noindent {\bf Remark 6.14}\,\,\,Combing (\ref{equa62}) with
(\ref{equa611}) we get an explicit formula for the number of simple
$\H_{q}(p,p,n)$-modules in the case when $\eps=q^{\ell}$ is a
primitive $p$-th root of unity, and $q$ is a primitive $p\ell$-th
root of unity. Note that (by Theorem \ref{main3}) this formula is
valid over any field $K$ as long as $K$ contains primitive $p$-th
root of unity and $\H_q(p,p,n)$ is split over $K$. Our formula
generalizes earlier results of Geck \cite{Ge} on the Hecke algebra
of type $D_n$ (i.e., of type $G(2,2,n)$). Note that Geck's method
depends on explicit information on character tables and
Kazhdan--Lusztig theory for Iwahori--Hecke algebras associated to
finite Weyl group, which are not presently available in our general
$G(p,p,n)$ cases. \medskip

We now deal with the remaining cases. That is, we assume the
following setup: $K=\mathbb{C}$, $q,\eps\in K$. $q$ is a primitive
$d\ell$-th root of unity, $q^{\ell}=\eps^{k}$ is a primitive $d$-th
root of unity, and $1\leq k<p$ is the smallest positive integer such
that $\eps^{k}\in\langle q\rangle$. In this case, we fix the order
of parameters $\{1,\eps,\cdots,\eps^{p-1}\}$ as
$\bbQ:=(\bbQ_1,\cdots,\bbQ_k)$, where
$\bbQ_i=(\veps^{i-1},\veps^{k+i-1},\cdots, \veps^{(d-1)k+i-1})$ for
$i=1,2,\cdots,k$. Let $\lam\in{K}_n$. We write
$\lam=(\lam^{[1]},\cdots,\lam^{[k]})$ (concatenation of ordered
tuples), where each $\lam^{[i]}$ is a $d$-multipartition. By Theorem
\ref{main2}, $$
\HH(\lam)=\Bigl(\HH'(\lam^{[k]}),\lam^{[1]},\cdots,\lam^{[k-1]}\Bigr),
$$
where $\HH'$ is as defined in Corollary \ref{maincor}.

As before, it suffices to derive a formula for $\widetilde{N}(m)$
for each integer $1\leq m<p$ satisfying $m|p$. Suppose that $1\leq
a\leq\min\{m,k\}$ is the greatest common divisor of $m$ and $k$.
By the division algorithm, we know that there exist integers $r_1,
r_2$ such that $a=r_1k+r_2m$. Let $\lam\in{K}_n$ be such
that $\HH^m(\lam)=\lam$.
Then $$\begin{aligned} &\quad
\Bigl((\HH')^{r_1}\bigl(\lam^{[1]}\bigr),
(\HH')^{r_1}\bigl(\lam^{[2]}\bigr),\cdots,(\HH')^{r_1}\bigl(\lam^{[k]}\bigr)\Bigr)
=\HH^{r_1k}(\lam)=\HH^{a-r_2m}(\lam)=\HH^a(\lam)\\
&=\Bigl((\HH')\bigl(\lam^{[k-a+1]}\bigr),(\HH')\bigl(\lam^{[k-a+2]}\bigr),\cdots,
(\HH')\bigl(\lam^{[k]}\bigr),\lam^{[1]},\lam^{[2]},\cdots,\lam^{[k-a]}\Bigr).\end{aligned}
$$
It follows that
\addtocounter{equation}{3}\begin{equation}\label{equa614}
(\HH')^{r_1}\bigl(\lam^{[i]}\bigr)=\begin{cases}
(\HH')\bigl(\lam^{[k-a+i]}\bigr), &\text{if $1\leq
i\leq a,$}\\
\lam^{[i-a]},&\text{if $a+1\leq i\leq k$.}\end{cases}
\end{equation}

By an easy induction, we get that $$
(\HH')^{-r_2m/a}\bigl(\lam^{[i]}\bigr)=(\HH')^{r_1k/a-1}\bigl(\lam^{[i]}\bigr)=\lam^{[i]},\,\,\forall\,1\leq
i\leq k.
$$

We claim that there exist positive integers $r'_1, r'_2$ such that
$a=r'_1k+r'_2m$ and $(r_2,r'_2)=1$. In fact, since $a|m$, it is easy
to check that
$a=\bigl(-(m/a-1)r_1\bigr)k+\bigl(-(m/a-1)r_2+1\bigr)m$. Taking
$r'_1=-(m/a-1)r_1$, $r'_2=-(m/a-1)r_2+1$, we prove our claim.

Now applying our previous argument again, we get that $$
(\HH')^{-r'_2m/a}\bigl(\lam^{[i]}\bigr)=\lam^{[i]},\,\,\forall\,1\leq
i\leq k.
$$
Let $x,y$ be two integers such that $xr_2+yr'_2=1$. Then $$
(\HH')^{-m/a}\bigl(\lam^{[i]}\bigr)=(\HH')^{-xr_2m/a-yr'_2m/a}\bigl(\lam^{[i]}\bigr)=
\lam^{[i]},\,\,\forall\,1\leq i\leq k.
$$
Equivalently,
$(\HH')^{m/a}\bigl(\lam^{[i]}\bigr)=\lam^{[i]},\,\,\forall\,1\leq
i\leq k$.
We write $\lam=(\lam^{[1]},\cdots,\lam^{[k]})$ (concatenation of ordered
tuples), where each $\lam^{[i]}$ is a $d$-multipartition. Let $n_i:=|\lam^{[i]}|$ for each $1\leq i\leq k$. Then $\sum_{i=1}^{k}n_i=n$. Moreover, in this case we see (from (\ref{equa614})) that $$
\lam^{[k-ja+i]}=(\HH')^{jr_1-1}\bigl(\lam^{[i]}\bigr),\,\,\text{for}\,\,
j=1,2,\cdots,k/a.
$$
As a consequence, $n_1+\cdots +n_a=na/k$.

By Theorem \ref{thm23}, Lemma \ref{lm24}, Theorem \ref{Ariki} and
our definition of the $p$-tuple $\bbQ$, we know that $\lam\in{K}_n$
if and only if for each $1\leq i\leq k$,
$\lam^{[i]}\in\mathcal{K}_{n_i}$, where $\mathcal{K}_{n_i}$ denotes
the set of Kleshchev $d$-multipartitions of $n_i$ with respect to
$(q,1,\eps',\cdots,{\eps'}^{d-1})$, $\eps'=\eps^k$. We define
$$\begin{aligned}
&\widetilde{\Sigma}(k,m):=\biggl\{\bigl(\lam^{[1]},\cdots,\lam^{[a]}\bigr)\vdash
\frac{na}{k}\biggm|\begin{matrix}
\text{$\lam^{[i]}\in\mathcal{K}_{n_i}, (\HH')^{m/a}({\lam^{[i]}})={\lam^{[i]}}$,}\\
\text{$\forall\,1\leq i\leq a,\,\,\sum_{i=1}^{a}n_i=\frac{na}{k}$}\end{matrix}\biggr\},\\
&\widetilde{N}(k,m):=\#\widetilde{\Sigma}(k,m).\end{aligned}
$$

\addtocounter{lem}{8}
\begin{lem} \label{lm616} With the above notations, the map which sends
$\lam=\bigl(\lam^{[1]},\cdots,$ $\lam^{[k]}\bigr)$ to
$\overline{\lam}:=\bigl(\lam^{[1]},\cdots,\lam^{[a]}\bigr)$
defines a bijection from the set $\widetilde{\Sigma}(m)$ onto the
set $\widetilde{\Sigma}(k,m)$.\end{lem}

\noindent {Proof:} \,Our previous discussion shows that the $p$-multipartition $\lam=(\lam^{[1]},\lam^{[2]},$ $\cdots,\lam^{[k]})\in\widetilde{\Sigma}(m)$ can be recovered from the $da$-multipartition $\overline{\lam}:=(\lam^{[1]},\lam^{[2]},$ $\cdots,\lam^{[a]})$ and the automorphism $\HH'$. In particular, the above map is injective.
It remains to prove the map is surjective.

Let $\alpha:=\bigl(\lam^{[1]},\cdots,\lam^{[a]}\bigr)\in
\widetilde{\Sigma}(k,m)$. Recall that $r_1k+r_2m=a$.
For each integer $1\leq i\leq a$, we define
$$
\lam^{[k-ja+i]}:=(\HH')^{jr_1-1}\bigl(\lam^{[i]}\bigr),\,\,\text{for}\,\,
j=1,2,\cdots,k/a.
$$
This is well-defined, since
$(\HH')^{kr_1/a-1}\bigl(\lam^{[i]}\bigr)=(\HH')^{-r_2m/a}\bigl(\lam^{[i]}\bigr)=\lam^{[i]}$
(because $\alpha\in\widetilde{\Sigma}(k,m)$). Note also
that the above definition is equivalent to (\ref{equa614}).
Therefore, if we set $\lam:=\bigl(\lam^{[1]},\cdots,$
$\lam^{[k]}\bigr)$, then the discussion above (\ref{equa614})
implies that $\HH^{r_2m}(\lam)=\lam$. Now recall that $$
r'_1=-(m/a-1)r_1,\,\, r'_2=-(m/a-1)r_2+1,\,\,a=r'_1k+r'_2m.
$$
Therefore, for each integer $1\leq i\leq a$, we have that
$$
(\HH')^{jr'_1-1}\bigl(\lam^{[i]}\bigr)=(\HH')^{-j(m/a-1)r_1-1}\bigl(\lam^{[i]}\bigr)=
\lam^{[k-ja+i]},\,\,\text{for}\,\,
j=1,2,\cdots,k/a.
$$
Therefore, the discussion above (\ref{equa614}) also implies that
$\HH^{r'_2m}(\lam)=\lam$. Since $xr_2+yr'_2=1$, it follows that
$\HH^m(\lam)=\HH^{xr_2m+yr'_2m}(\lam)=\lam$. In other words,
$\lam\in\widetilde{\Sigma}(m)$ with $\overline{\lam}=\alpha$, as required. This proves the
surjectivity, and hence completes the proof of the lemma.\hfill\qed
\medskip

Let $\widetilde{d}:=(d,\frac{m}{a})$. Note that $(\HH')^{d}=\id$.
Therefore, $(\HH')^{m/a}(\lam^{[i]})=\lam^{[i]}$ if and only
$(\HH')^{\widetilde{d}}(\lam^{[i]})=\lam^{[i]}$. Now applying
Lemma \ref{lm616}, Theorem \ref{mainthm3} and (\ref{equa63}), we
get that

\addtocounter{thm}{6}
\begin{thm} \label{mainthm4} Suppose that $K=\mathbb{C}$,\, $q,\eps\in K$. $q$ is a primitive
$d\ell$-th root of unity, $q^{\ell}=\eps^{k}$ is a primitive
$d$-th root of unity, and $1\leq k<p$ is the smallest positive
integer such that $\eps^{k}\in\langle q\rangle$. Let $1\leq m<p$
be an integer such that $m|p$. Let $a=(m,k)$, $\widetilde{d}:=(d,\frac{m}{a})$. Then
$$\begin{aligned}
&\quad\,\widetilde{N}(m)=\widetilde{N}(k,m)\\
&=\sum_{n_1+\cdots+n_{a}=\frac{na}{k}}\prod_{i=1}^{a}
\biggl(\#\Irr\H_{q'',1,\eps'',\cdots,(\eps'')^{\widetilde{d}-1}}\Bigl(\widetilde{d},\frac{\widetilde{d}n_i}{d}\Bigr)
\biggr),
\end{aligned}$$ where $\eps''=(q'')^{\l}$ is a primitive
$\widetilde{d}$-th root of unity, and $q''$ is a primitive $(\widetilde{d}\l)$-th
root of unity. In particular, for each integer $1\leq
\widetilde{m}<p$ such that $\widetilde{m}|p$, we get that
\addtocounter{equation}{2}
\begin{equation}\label{equa617}\begin{split}
&\quad\,N(\widetilde{m})=\sum_{1\leq m\leq\widetilde{m},
m|\widetilde{m}}\mu(\widetilde{m}/m)\widetilde{N}(m)\\
&=\sum_{1\leq m\leq\widetilde{m},
m|\widetilde{m}}\mu(\widetilde{m}/m)
\sum_{n_1+\cdots+n_{a}=\frac{na}{k}}\prod_{i=1}^{a}
\biggl(\#\Irr\H_{q'',1,\eps'',\cdots,(\eps'')^{\widetilde{d}-1}}\Bigl(\widetilde{d},\frac{\widetilde{d}n_i}{d}\Bigr)
\biggr).\end{split}\end{equation}
\end{thm}

\noindent {\bf Remark 6.19}\,\,\,Combing (\ref{equa62}) with
(\ref{equa617}) we get an explicit formula for the number of simple
$\H_{q}(p,p,n)$-modules in the case when $q$ is a primitive
$d\ell$-th root of unity, $q^{\ell}=\eps^{k}$ is a primitive $d$-th
root of unity, and $1\leq k<p$ is the smallest positive integer such
that $\eps^{k}\in\langle q\rangle$. Thus the problem on determining
explicit formula for the number of simple $\H_{q}(p,p,n)$-modules is
solved by our Theorem \ref{mainthm3} and Theorem \ref{mainthm4} in
all cases. As before, this formula is valid over any field $K$ as
long as $K$ contains primitive $p$-th root of unity and
$\H_q(p,p,n)$ is split over $K$. Finally, we remark that, these
explicit formulas strongly indicate that there are some new intimate
connections between the representation of $\H_{q}(p,p,n)$ at roots
of unity and the representation of various Ariki--Koike algebras of
smaller sizes at various roots of unity. It seems very likely that
the decomposition matrix of the latter can be naturally embedded as
a submatrix of the decomposition matrix of the former. We will leave
it to a future project.

\end{document}